\documentclass[a4paper]{lipics-v2021}

\nolinenumbers
\title{{A}cyclic {C}oloring of {P}roducts of {D}igraphs and of {D}igraphs with {B}ounded {T}reewidth}

\author{Isnard Lopes Costa}{ParGO Group - Parallelism, Graphs and Optimization \and Departamento de Matem\'atica - Universidade Federal do Cear\'a, Fortaleza, Brazil}{isnard.lopes@alu.ufc.br}{https://orcid.org/0000-0002-9478-2378}{}

\author{Ana Shirley Silva }{ParGO Group - Parallelism, Graphs and Optimization \and Departamento de Matem\'atica - Universidade Federal do Cear\'a, Fortaleza, Brazil}{anasilva@mat.ufc.br}{https://orcid.org/0000-0001-8917-0564}{}

\authorrunning{I. Costa and A. S. Silva} 

\Copyright{Isnard Costa and Ana Silva} 

\ccsdesc[100]{Mathematics of computing~Graph theory}
\ccsdesc[100]{Mathematics of computing~Graph coloring} 

\keywords{dichromatic number, products of graphs, products of cycles, treewidth, parameterized complexity.} 

\funding{Partially supported by CAPES and CNPq, Brazil}

\EventEditors{John Q. Open and Joan R. Access}
\EventNoEds{2}
\EventLongTitle{42nd Conference on Very Important Topics (CVIT 2016)}
\EventShortTitle{CVIT 2016}
\EventAcronym{CVIT}
\EventYear{2016}
\EventDate{December 24--27, 2016}
\EventLocation{Little Whinging, United Kingdom}
\EventLogo{}
\SeriesVolume{42}
\ArticleNo{23}

\usepackage{graphicx}
\usepackage{amsmath,amsthm}
\usepackage{float}
\usepackage{thmtools}
\usepackage{thm-restate}
\usepackage{pdfpages}
\usepackage{hyperref}
\usepackage{ amssymb }
\usepackage{pgf}
\usepackage{caption}
\usepackage{tikz}
\usepackage{graphicx}
\usetikzlibrary{arrows}
\usepackage{complexity}

\newtheorem{question}{Question}
\newif\iffigures
    \figurestrue 

\newcommand{\cn}{\chi_a}
\newcommand{\vv}{\overset{\rightarrow}}
\newcommand{\vvv}{\overset{\leftrightarrow}}

\begin{document}
\maketitle

\begin{abstract}

  The dichromatic number of a digraph $G$ is the smallest integer $\cn(G)$ such that the vertex set of $G$ can be partitioned into $\cn(G)$ sets, each of which induces an acyclic subdigraph. This is a generalization of the classic chromatic number of graphs. Here, we investigate the dichromatic number of the cartesian, direct, strong and lexicographic products, giving generalizations of some classic results on the chromatic number of products. More specifically, we prove that the following inequalities, known to hold for the chromatic number of graphs, still hold for the dichromatic number of digraphs:  $\cn(G\square H)=\max\{\cn(G),\cn(H)\}$;  $\cn(G\times H)\le \min\{\cn(G),\cn(H)\}$; and $\cn(G[H]) = \cn(G[\vvv{K}_k])$, where $k =\cn(H)$ and $\vvv{K}_k$ denotes the complete digraph on $k$ vertices. In addition, we investigate the products of directed cycles, giving exact values for  $\cn(\vv{C}_n\times \vv{C}_m)$ and $\cn(\vv{C}_n\boxtimes \vv{C}_m)$ for every $n,m$, and for $\cn(\vv{C}_n[H])$ for every positive integer $n$. This latter result generalizes a result given in~\cite{PP.16}, where they give exact values when $n>\cn(H)$. We also provide a upper-bound to the dichromatic number of a digraph $G$ as a function of the treewidth of its underlying graph and we present an {\FPT}-time algorithm that computes the dichromatic number of $G$, when parameterized by treewidth of the underlying graph of $G$.

\end{abstract}

\section{Introduction}\label{intro}

In this work, we consider only digraphs with no loops and no multiple arcs. We also make many observations about colorings of simple undirected graphs, which we call simply \emph{graphs}~\footnote{We refer the reader to~\cite{bang2008digraphs} for basic definitions.}. 
Given a graph $G$, a \emph{proper coloring of $G$} is a partition of its vertex set into independent sets (sets containing no edges), and the \emph{chromatic number of $G$} is the minimum value $k$ for which $G$ admits a proper coloring with $k$ colors; it is denoted by $\chi(G)$. 
Given a digraph $D$, an \emph{acyclic coloring} of $D$ is a partition of its vertex set into  subsets that induce acyclic subdigraphs of $D$. The smallest integer $k$ such that $D$ has an acyclic coloring that uses $k$ colors is called the \emph{dichromatic number} of $D$ and is denoted by $\cn(D)$. Observe that, given a graph $G=(V,E)$, we can obtain a digraph  $D=(V,E')$ related to $G$ such that, for each edge $uv\in E$, there are two arcs $uv$ and $vu$ in $E'$; this is called the \emph{symmetric digraph of $G$} and is denoted by $D(G)$. Similarly, we say that $G$ is the \emph{underlying graph of $D$} and denote it $U(D)$. Because each edge of $G$ gives rise to a cycle of length $2$ in $D(G)$, we get that an acyclic coloring of $D(G)$ is equivalent to a proper coloring of $G$. This means that the concept of acyclic coloring in digraphs generalizes the concept of proper colorings in  graphs and therefore inherits all of the {$\NP$}-completeness results. For instance, it follows that it is $\NP$-complete to decide, given a digraph $D$, whether $\cn(D)\le 3$~\cite{H.81}. However, it is also known that deciding $\cn(D)\le 2$ is $\NP$-complete~\cite{BFJKM.04}, even if $D$ is an \emph{oriented graph} (digraph with no cycles of length~2), unlike the classical coloring problem on graphs, where it is polynomial to decide $\chi(G)\le 2$ for every $G$.  Acyclic coloring was introduced in \cite{NL.82}, and has attracted many interest in the past years  (see e.g.~\cite{ACH.16, B.etal.13, FHM.03}).

Because these colorings are generalizations of the classical proper colorings, many authors are interested in adapting classical results about the chromatic number of graphs. For example, it is well known that the chromatic number of a graph $G$ is bounded above by its maximum degree plus~1. 
This has been generalized to acyclic colorings in \cite{NL.82}, where the author proves that the dichromatic number of a digraph is at most 1 plus the minimum between the maximum outdegree and the maximum indegree of a vertex. 
In this work, we continue this line of investigation, presenting the acyclic coloring counterpart of classical results about proper colorings of products of graphs. We also give an algorithm to compute the dicromatic number that is $\FPT$ when parameterized by the treewith of the underlying graph. We mention that $\FPT$ algorithms (and sometimes even $\XP$ algorithms) parameterized by other directed parameters (as for instance directed treewidth, DAG width and directed cliquewidth) are not expected to exist according to results in~\cite{GKR.20,MSW.19}. We mention additionally that these results have been presented in the Msc. Thesis of Isnard Lopes~\cite{L.20}, and that, since then, a better FPT algorithm parameterized by treewidth has been proposed by Harutyunyan, Lampis and Melissinos~\cite{HLM.21}.

Our paper is organized as follows. In Section~\ref{sec:prelim}, we give basic concepts and present our results regarding products of graphs, making a parallel with existing results on proper colorings. In Sections~\ref{sec:general} and~\ref{sec:dicycles}, we present the proofs of said results. Finally, in Section~\ref{sec:treewidth} we present our $\FPT$ algorithm and in Section~\ref{sec:conclusion} we make our final remarks.

\section{Preliminaries and Related Work}\label{sec:prelim}

A product between digraphs $G$ and $H$ is a digraph defined over the set of vertices $ V (G) \times V (H) $. Below, we define the arc sets of the digraph products studied here (see Figures~\ref{fig:d3prod} and ~\ref{fig:dlexic} for illustrations). We call digraphs $G$ and $H$ the \emph{factors} of the product.

	\begin{enumerate}
	\item Cartesian product: $(u,x)(v,y)\in E(G~\square~H)$ if, and only if, $u=v$ and $xy \in E(H)$, or $uv \in E(G)$ and $x=y$;
    \item  Direct product: $(u,x)(v,y)\in E(G\times H)$ if, and only if,  $uv \in E(G)$ and $xy \in E(H)$;
    \item  Strong product: $(u,x)(v,y)\in E(G\boxtimes H)$ if, and only if, $u=v$ and $xy \in E(H)$, or $uv \in E(G)$ and $x=y$, or $uv \in E(G)$ and $xy \in E(H)$;
    \item  Lexicographic product: $(u,x)(v,y)\in E(G[H])$ if, and only if, $u=v$ and $xy \in E(H)$, or $uv \in E(G)$.
	\end{enumerate}

\iffigures
	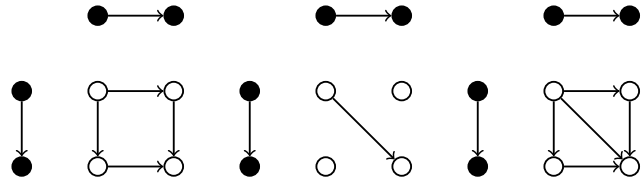
\begin{figure}[H]
	    \centering

	\begin{tikzpicture}

 \pgfsetlinewidth{.7pt}

  \tikzstyle{blackvertex}=[draw,circle,fill=black,minimum size=7pt,inner sep=0pt]
   \tikzstyle{vertex}=[draw,circle,fill=white,minimum size=7pt,inner sep=0pt]
   \tikzstyle{grayvertex}=[draw,circle,fill=black!25,minimum size=7pt,inner sep=0pt]

	\pgfsetarrows{->}
	\node[blackvertex] (1) at (0,1){};
	\node[blackvertex] (2) at (0,0){};
	\draw[->] (1) -> (2);
	
	\node[blackvertex] (3) at (2,2){};
	\node[blackvertex] (4) at (1,2){}; 
	\draw[->] (4)--(3);	
	
	\node[vertex](5) at (1,1){};
	\node[vertex](6) at (2,1){};
	\node[vertex](7) at (2,0){};
	\node[vertex](8) at (1,0){};

	\node[blackvertex] (11) at (3,1){};
	\node[blackvertex] (22)at (3,0){};
	\draw[->]  (11)--(22);
	
	\node[blackvertex] (33)at (5,2){};
	\node[blackvertex] (44) at (4,2){}; 
	\draw (44)--(33);		
	
	\draw[->] (5) -- (6);
	\draw[->] (6) -- (7);
	\draw[->] (8) -- (7);
	\draw[->] (5) -- (8);
	
	\node[blackvertex] (111) at (6,1){};
	\node[blackvertex] (222)at (6,0){};
	\draw[->]  (111)--(222);
	
	\node[blackvertex] (333)at (7,2){};
	\node[blackvertex] (444) at (8,2){}; 
	\draw[->] (333)--(444);

    \node[vertex](9) at (4,1){};
	\node[vertex](10) at (5,1){};
	\node[vertex](11) at (5,0){};
	\node[vertex](12) at (4,0){};
	
	\draw[->] (9)--(11);

    \node[vertex](13) at (7,1){};
	\node[vertex](14) at (8,1){};
	\node[vertex](15) at (8,0){};
	\node[vertex](16) at (7,0){};
	
	\draw[->] (13) -- (14);
	\draw[->] (14) -- (15);	
	\draw[->] (16) -- (15);
    \draw[->] (13) -- (15);
	\draw[->] (13) -- (16);
	
	\end{tikzpicture}
    \vspace{.5cm}

	    \caption{This figure illustrates the cartesian, direct and strong products between two dipaths $\vv{P}_2 $'s. The vertices of the factors are highlighted in black. From left to right, we have $ \vv{P_2} ~ \square ~ \vv{P_2} $, $ \vv{P_2} \times \vv{P_2} $, $ \vv{P_2} \boxtimes \vv{P_2} $.}
	    \label{fig:d3prod}
	\end{figure}
	
	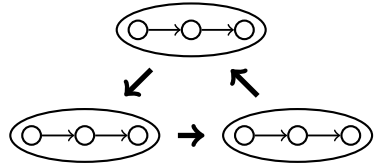
\begin{figure}[H]
	    \centering

 \pgfsetlinewidth{1pt}

  \tikzstyle{blackvertex}=[draw,circle,fill=black,minimum size=7pt,inner sep=0pt]
   \tikzstyle{vertex}=[draw,thick,circle,fill=white,minimum size=7pt,inner sep=0pt]
   \tikzstyle{grayvertex}=[draw,circle,fill=black!25,minimum size=7pt,inner sep=0pt]

\begin{tikzpicture}[scale=.7]

 \pgfsetlinewidth{.7pt}

  \tikzstyle{blackvertex}=[draw,circle,fill=black,minimum size=7pt,inner sep=0pt]
   \tikzstyle{vertex}=[draw,thick,circle,fill=white,minimum size=7pt,inner sep=0pt]

             \node[vertex](ax) at (1,2){};
             \node[vertex](ay) at (2,2){};
             \node[vertex](az) at (3,2){};
             
             \draw[thick] (0,0) ellipse (1.4cm and .5cm);
             
             \node[vertex](bx) at (-1,0){};
             \node[vertex](by) at (0,0){};
             \node[vertex](bz) at (1,0){};
             
             \draw[thick] (4,0) ellipse (1.4cm and .5cm);
             
             \node[vertex](cx) at (3,0){};
             \node[vertex](cy) at (4,0){};
             \node[vertex](cz) at (5,0){};
             
             \draw[thick] (2,2) ellipse (1.4cm and .5cm);

             \draw[->] (ax) -- (ay);
             \draw[->] (ay) -- (az);
             \draw[->] (cx) -- (cy);
             \draw[->] (cy) -- (cz);
             \draw[->] (bx) -- (by);
             \draw[->] (by) -- (bz);
             
             {
             \pgfsetlinewidth{2pt}
             
             \draw[->] (1.25,1.25) -- (0.75,0.75);
             \draw[->] (1.75,0) -- (2.25,0);
             \draw[->] (3.25,0.75) -- (2.75,1.25);
             }

\end{tikzpicture}

	    \caption{Illustration of the digraph $\vv{C}_3[\vv{P}_3]$. The thick arcs represent the existence of all possible arcs between the corresponding sets.}
	    \label{fig:dlexic}
	\end{figure}
\fi	
	
These products are defined in the same way for undirected graphs, and we mention that the symbols $\square$, $\times$ and $\boxtimes$ are useful to remember how the product looks like, since these symbols represent also the drawings that one gets when picking the related products on a pair of edges of undirected graphs. One should also observe that, if $G$ and $H$ are graphs, then $D(G\odot H) = D(G)\odot D(H)$ for each possible product $\odot$. Therefore, results that hold for the acyclic coloring of products of general digraphs also hold for proper colorings of products of graphs. 

In what follows, we discuss our results obtained on each of these products. We denote a directed path (also called \emph{dipath}) on $n$ vertices by $\vv{P}_n$, a directed cycle (also called \emph{dicycle}) on $n$ vertices by $\vv{C}_n$, and the complete digraph on $n$ vertices by $\vvv{K}_n$. Additionally, given a graph $G$, an \emph{orientation of $G$} is any digraph $D$ such that exactly one between $uv$ and $vu$ is in $A(D)$ for every $uv\in E(G)$.

We start by citing the following result, which holds for the chromatic number of graphs and has been proved to also hold for acyclic colorings of digraphs.

\begin{theorem}[\cite{PP.16}]\label{prop:produto} Let $G$ and $H$ be digraphs. Then, $\cn(G[H]) \leq \cn(G)\cdot \cn (H)$.
\end{theorem}

The interesting aspect of the above result is that it helps us solve the problem whenever one of the factors is a DAG (directed acyclic graph). Indeed, it is known that the lexicographic product contains all the other products, and that the lexicographic, cartesian and strong products contain copies of both factors. Also, observe that the dichromatic number is a monotonic metric, i.e., $\cn(H)\le \cn(G)$ whenever $H\subseteq G$. We then get that, if $H$ is a DAG, then the dichromatic number of $G$ times $H$ equals the one of $G$ for each of the aforementioned products. Additionally, by Theorem~\ref{thm:folkore} presented later, we get that $\cn(G\times H) = 1$ whenever $H$ is a DAG.  
Therefore, the problem is not so interesting when we restrict ourselves to the case where one of the factors is a DAG. In particular, this includes orientations of paths, or orientations of cycles that are not a dicycle. This is why in what follows we are interested only in dicycles, and not in acyclic orientations of cycles. Below, we formally state the mentioned result.

\begin{corollary}\label{cor:DAG} 
Let $G$ and $H$ be digraphs, and $\odot$ be a digraph product. If either $G$ or $H$ is a DAG, then: $\cn(G\odot H) = \cn(G)$ if $\odot$ is the cartesian, strong or lexicographic product; otherwise, we get $\cn(G\times H) = 1$. In particular, this holds when either $G$ or $H$ is an oriented path, or an oriented cycle that is not a dicycle. 
\end{corollary}

A classical result for the chromatic number of the cartesian product of two graphs $G$ and $H$ says that it is the maximum between $\chi(G)$ and $\chi(H)$~\cite {S.57}. We show here that the same is true when it comes to the product of two digraphs. 

\begin{restatable}{theorem}{thmCartesian}\label{thm:Cartesian}
Let $G$ and $H$ be digraphs. Then, \[\cn(G~\square~H) = \max\{\cn(G),\cn(H)\}.\]
\end{restatable}

In contrast, the minimum over $\chi(G)$ and $\chi(H)$ is known to be an upper bound for the chromatic number of $G\times H$ (folklore). Here, we show that this also holds for the dichromatic number:
		
	\begin{restatable}{theorem}{thmfolklore}\label{thm:folkore} Let $G$ and $H$ be digraphs. Then, \[\cn(G \times H) \le \min \{\cn(G), \cn(H)\}.\]
	
	\end{restatable}

A famous conjecture about the chromatic number of the direct product of graphs is the {H}edetniemi's {C}onjecture~\cite {H.66}, which stated that the above bound is tight, i.e., that $\chi(G \times H) = \min \{\chi(G),\chi(H)\}$. This holds for graphs with chromatic number at most $4$~\cite{ZS.85}, but the conjecture remained open for 52 years until it was recently proved false by {S}hitov~\cite{S.19}. Here, we present a result similar to the one in~\cite{ZS.85}, proving that equality holds for digraphs with dichromatic number at most~2.

\begin{restatable}{theorem}{thmpositiveAnswer}\label{thm:positiveAnswer}
Let $G,H$ be digraphs. If $\min\{\cn(G),\cn(H)\} = \cn(G) \le 2$, then $\cn(G\times H) = \cn(G)$.
\end{restatable}

Regarding the strong product, up to our knowledge, there are no general upper bounds. Nevertheless, 
the following results study the chromatic number of the strong product of cycles, complete graphs and some {K}neser graphs, giving exact results for certain values: 

\begin{theorem}
    \begin{enumerate}
        \item (\cite{HR.82},\cite{V.79}) For $k \ge 2 $ and $n \ge 2$, $\chi(C_{2k+1} \boxtimes C_{2n+1})=5$.
        \item (\cite{K.93}) For $k \ge 2$, $\chi(C_{5} \boxtimes C_5 \boxtimes C_{2k+1}) = 10 + \lceil \frac{5}{k} \rceil$.
        \item (\cite{S.76}) For $k\ge 2$ and $n \ge 1$, $\chi(C_{2k+1} \boxtimes K_n )= 2n + \lceil \frac{n}{k}\rceil$. And if $H$ is any bipartite graph, then $\chi(H \boxtimes K_n )= 2n$
        \item (\cite{KM.94}) For $k\ge 0$ and $n \ge 1$, $\chi (KN_{n,k} \boxtimes K_n)= 2n + k $
    \end{enumerate}
\end{theorem}
    
Complementing the above results for the chromatic number of graphs, here we give an exact value for the product $H\boxtimes C_{2k+1}$ for when $H$ is a bipartite graph. In~\cite{V.79}, the author proves that the chromatic number of the strong product is at least the maximum between the chromatic numbers of the factors plus $2$. This implies that if both factors are even cycles (or bipartite, in general), then the chromatic number of the strong product is at least $4$. And since $\chi(G\boxtimes H)\le \chi(G[H])\le \chi(G)\cdot\chi(H)$, we get that $\chi(C_{2k}\boxtimes C_{2\ell}) = 4$. Observe that combining these results, we get the exact value for $\chi(C_m\boxtimes C_n)$ for every choice of $m$ and $n$.

\begin{restatable}{theorem}{thmprodforteoutros}\label{thm:ProdForte_outros}
  Let $n \ge 2$ be an integer, and $H$ be a bipartite graph with at least one edge. Then, $\chi(H \boxtimes C_{3})=6$ and $\chi(H \boxtimes C_{2n+1})=5$. 
\end{restatable}

We also give exact values for the dichromatic number of the strong product of two dicycles.

\begin{restatable}{theorem}{thmprodforte}\label{thm:prodforte} 
         Let $\vv{C}_m, \vv{C}_n$ be dicycles with $m, n$ vertices, respectively. Then,
         \begin{enumerate} 
         \item If $m=n=2$, then $\cn(\vv{C}_m \boxtimes \vv{C}_n) = 4$,
         \item If $m = 3$ and $ n \in\{2, 3\}$, then $\cn(\vv{C}_m \boxtimes \vv{C}_n) = 3$,
         \item If $n \ge 3$ and  $ m \ge 4$, then $ \cn(\vv{C}_m \boxtimes \vv{C}_n) = 2$.
         \end{enumerate}
\end{restatable}

\medskip

Regarding the lexicographic product, in \cite{GS.75} the authors proved the equivalence between coloring the lexicographic product of $G$ by $H$ and coloring the product of $G$ by the complete graph with $\chi(H)$ vertices.
We show that the analogous result holds also for  acyclic colorings.

\begin{restatable}{theorem}{thmGbyKk}\label{thm:GbyKk}
         Let $G,H$ be digraphs, with $\cn(H)=k$. Then $\cn(G[H]) = \cn(G[\vvv{K}_k]).$
\end{restatable}

Still concerning the lexicographic product, in the seminal work \cite{NL.82}, among other results, the author establishes that the dichromatic number of $ G[H] $ is at least $ \cn (G) + \cn (H) -1 $. In \cite{PP.16}, the authors prove that this bound is tight when $G$ is the dicycle with $n$ vertices and $H$ is an arbitrary digraph with $n> \cn(H)$; that is, they prove that if $n> \cn(H)$, then $\cn (\vv{C}_n [H]) = \cn(H) + 1 $. We generalize their result, giving an exact value for every $n$.

\begin{restatable}{theorem}{thmnlechi}\label{thm:nlechi}
Let $n$ be any positive integer and $\vv{C}_n$ be a dicycle on $n$ vertices. For every digraph $H$ we have that:
         \[\cn(\vv{C}_n[H]) = \cn(H) + \left\lceil  \frac{\cn(H)}{n-1} \right\rceil.\]
\end{restatable}

Before we proceed to our proofs, we need some further definitions, notation and basic results.

Here, we denote by $[k]$ the set of the integer $\{1,\ldots, k\}$. Also, given a function $f:X \to Y$ and $X'\subseteq X$, we denote by $f(X')$ the set $\{f(x)\mid x\in X'\}$, and for $y\in Y$, we denote by $f^{-1}(y)$ the set $\{x \in X \mid f(x)=y\}$.

Given a digraph $G$, a subset $S\subseteq V(G)$ is called \emph{acyclic} if the subdigraph of $G$ induced by $S$ is an acyclic digraph (also called DAG). A subdigraph $H$ of $G$ is \emph{strong} if for every $u,v \in V(H)$ there exists a dipath from $u$ to $v$ in $H$, and a dipath from $v$ to $u$ in $H$; and it is a \emph{strong component of $G$} if it is a maximal strong subdigraph of $G$. 

Given digraphs $G$ and $H$, an \emph{homomorphism} of $ G $ into $ H $ is a function that preserves the arcs; more formally, it is a function $ f: V (G) \rightarrow V (H) $ such that $ f(u)f (v) \in E (H) $, for every $ uv \in E (G) $. Note that when $ H $ is the complete digraph $\vvv{K}_k $, then an homomorphism of $D(G)$ into $H$ is also a proper $k$-coloring of $G$, and vice versa; that is, the notion of homomorphism generalizes proper colorings. However, if we want a generalization of acyclic colorings we need an additional constraint. 
An \emph{acyclic homomorphism} of a digraph $G$ into a digraph $H$ is an homomorphism $f$ such that $f^{-1}(x)$ is acyclic in $G$ for every $x\in V(H)$. It follows by definition that $\cn(G)\le k$ if and only if $G$ has an acyclic homomorphism into $\vvv{K}_k$. This notion was introduced in~\cite{FHM.03}, where the authors also prove the following proposition, which will be used later.

\begin{proposition}[\cite{FHM.03}]\label{prop:transitivity}
Let $G,H,F$ be digraphs. If $G$ has an acyclic homomorphism into $H$, and $H$ has an acyclic homomorphism into $F$, then $G$ has an acyclic homomorphism into $F$.
\end{proposition}

\medskip
The following lemmas will be useful in the remainder of the text.

\begin{lemma}\label{lemma:cycleinColumn}
Let $ G $, $H $ be digraphs, and $ C = (w_1, \ldots, w_t) $ be a dicycle in $G\boxtimes H$. Also, for each $ i \in [t] $, let $ w_i = (u_i, x_i) $, and let $U = \{u_1, \cdots, u_t\}$. 
If $|U|>1$, then $U$ contains a dicycle in $G$. The same holds if $C$ is a dicycle in $G\square H$.

	\end{lemma}
	\begin{proof}
	  
	Observe that, by the definition of strong product, for each $i\in [t-1]$, since $w_iw_{i+1}\in E(G\boxtimes H)$, we get that either $u_i=u_{i+1}$ or $u_iu_{i+1}\in E(H)$. Because this also holds for $u_t$ and $u_1$, and since $|U|>1$, we get that  $U$ contains some dicycle. The same argument can be applied for the cartesian product.
	\end{proof}
	
	\begin{lemma}\label{lemma:key}
Let $ G $ and $ H $ be digraphs, and $ C = (w_1, \ldots, w_t) $ be a dicycle in $G[H]$. Also, for each $ i \in [t] $, let $ w_i = (u_i, x_i) $. Then, either $U = \{u_1, \cdots, u_t\}$ contains a dicycle in $G$, or $|U|=1$ and $X = \{x_1,\cdots, x_t\}$ forms a dicycle in $H$.
	\end{lemma}
	
		\begin{proof}
	First, note that if $U=\{u\}$, then, by the definition of lexicographic product, we have that  $x_ix_{(i+1) \mod t} \in E(H)$, for every $i \in [t]$, i.e., $X$ forms a dicycle in $H$. 
	Now,  suppose that $|U| \ge 2$. Let $u_{t_1},\ldots,u_{t_r}$, with $r\in [t]$, be the vertices of $U$ disregarding consecutive repetitions, and enumerate them in the order of appearance in $C$ (i.e. $t_1 < t_2 < \cdots t_r$). 
	By the definition of lexicographic product, since $w_{t_i}w_{t_{(i+1) \mod r}} \in E(G[H])$ and $u_{t_i} \neq u_{t_{(i+1) \mod r}}$, we have $u_{t_i}u_{t_{(i+1) \mod r}} \in E(G)$, for each $i \in [r]$. Hence, the subdigraph of $G$ induced by the vertices $u_{t_1},\ldots,u_{t_r}$ contains a dicycle. 
	\end{proof}

We first present the proofs of the general results (Section~\ref{sec:general}), and then we prove the results concerning dicycles (Section~\ref{sec:dicycles}).

\section{General Bounds for Digraph Products}\label{sec:general}

In this section, we prove the results concerning general bounds; more especifically, we prove Theorems~\ref{thm:Cartesian},~\ref{thm:folkore} and~\ref{thm:GbyKk}. We restate each theorem for clarity purposes.

  \thmCartesian*
\begin{proof}
Suppose, without loss of generality, that $\max\{\cn(G), \cn(H)\} = \cn(G) = k$. Since $ G \subseteq G~\square ~H$, we have $\cn(G~\square ~ H) \geq k$. Consider now acyclic colorings $f_1 : V(G) \rightarrow [k]$ and $f_2 : V(H) \rightarrow [k]$ of $G$ and $ H$, respectively. Define the coloring $f : V(G)\times V(H) \rightarrow \{ 0, \ldots , k-1 \}$ such that $f(u, x)= (f_1(u) + f_2(x)) \mod k$. We show that $f$ is an acyclic coloring of $G ~\square~H$, but before we proceed note that the following holds.
\begin{itemize}
    \item[(*)] For every $u\in V(G)$, and every $x,y\in V(H)$, we have $f(u,x) = f(u,y)$ if and only if $f_2(x) = f_2(y)$. Similarly, For every $x\in V(H)$, and every $u,v\in V(G)$, we have $f(u,x) = f(v,x)$ if and only if $f_1(u) = f_2(v)$.
\end{itemize}

Suppose by contradiction that $f$ yields a monochromatic dicycle $C=(w_1, \ldots , w_t)$ of color $c$ and write $w_i = (u_i, x_i)$, for each $i \in [t]$. 
Consider $U=\{u_1,\ldots, u_t\}$ and $X = \{x_1,\ldots,x_t\}$. We prove that $f_1(u_i) = f_1(u_j)$ and $f_2(x_i) = f_2(x_j)$, for every $i,j\in [t]$, thus getting a contradiction since $f_1$ and $f_2$ are acyclic colorings and, by Lemma \ref{lemma:cycleinColumn}, either $U$ or $X$ contains a dicycle. Let $d_1 = f_1(u_1)$, and $d_2 = f_2(x_1)$. By the definition of lexicographic product, and since $(u_1,x_1)(u_2,x_2)\in E(G~\square~H)$, we know that either $u_2 = u_1$, or $x_1=x_2$. If the former occurs, then clearly $f_1(u_2) = f_1(u_1) = d_1$ and, because $f(u_1,x_1) = f(u_1,x_2) = c$, it follows from (*) that $f_2(x_2) = f_2(x_1) = d_2$. The same can be obtained when $x_1 = x_2$. By iteratively applying this argument to the remaining vertices of $C$, since they form a dicycle in $G~\square~H$, it follows that $f_1(u_i) = d_1$ and $f_2(x_i) = d_2$ for every $i\in [t]$, as desired. 

\end{proof}

\thmfolklore*
	 \begin{proof}
 Let $f:V(G) \rightarrow [k]$ be an acyclic coloring of $G$, and define $g: V(G) \times V(H) \rightarrow [k]$ by $g(u, x)=f(u)$, for every $(u,x)\in V(G)\times V(H)$. We show that $g$ is an acyclic coloring of $G \times H$. Suppose by contradiction that some color class, say $c$, is not acyclic, and let  $(w_1, \ldots w_t)$ be a dicycle of color $c$ in $G \times H$; write $w_j=(u_j, x_j)$ for each $j \in [t]$. Note that, by the definition of direct product, we have that $u_iu_{(i+1)\mod t}\in E(G)$ and $x_ix_{(i+1)\mod t}\in E(H)$ for every $i\in [t]$. It thus follows that both $\{u_1,\cdots,u_t\}$ and $\{x_1,\cdots,x_t\}$ contain dicycles. We get a contradiction because $f(u_j)=g(w_j)=c=g(w_{j'})=f(u_{j'})$ for all $j, j' \in [t]$, i.e., $\{u_1, \ldots, u_t\}$ contains a monochromatic dicycle of color $c$ in $G$. 
 \end{proof}

In the following result we prove the equivalence between coloring the lexicographic product of $G$ by $H$ and coloring the lexicographic product of $G$ by the complete digraph with $\cn(H)$ vertices.

\thmGbyKk*

\begin{proof}
We proceed by showing that the two inequalities hold simultaneously.
\begin{itemize}

\item[($\leq$)] Since $ k= \cn(H)$, there exists an acyclic homomorphism $ f : V(H) \rightarrow V(\vvv{K}_k)$. Define $g: V(G[H]) \rightarrow V(G[\vvv{K}_k])$  such that $g(u, x)=(u, f(x))$. We show that $g$ is an acyclic homomorphism, thus getting $\cn(G[H])\le \cn(G[\vvv{K}_k])$, by Proposition~\ref{prop:transitivity}. Recall that $g$ is an acyclic homomorphism if and only if: for every $(u,x)(v,y) \in E(G[H])$ we have either $g(u,x)=g(v,y)$ or $g(u,x)g(v,y) \in E(G[\vvv{K}_k])$; and for each  $(u,i) \in V(G[\vvv{K}_k])$, the set $g^{-1}(u,i)$ is an acyclic set in $G[H]$. 
Consider an edge $(u, x)(v, y) \in E(G[H])$, and first suppose that $u=v$. By the definition of lexicographic product, we get that $xy \in E(H)$.  Furthermore, since $f$ is an acyclic homomorphism of $H$ into $\vvv{K}_k$, we have that either $f(x)=f(y)$ or $f(x)f(y) \in E(\vvv{K}_k)$. If the former occurs, we have $g(u, x)= (v, f(x)) = (u, f(y))= g(u, y)$. And if the latter occurs, we have $g(u, x)g(u, y) = (u, f(x))(u, f(y))$, which in turn is an edge in $G[\vvv{K}_k]$. Now, suppose $u \ne v$. By the definition of lexicographic product, we have $uv \in E(G)$, which in turn implies that $(u,i)(v,j) \in E(G[\vvv{K}_k])$, for every $i,j \in V(\vvv{K}_k)$. In particular, $(u, f(x))(v, f(y)) = g(u, x)g(v, y)  \in E(G[\vvv{K}_k])$.

Finally, let $(u, i) \in V(G) \times V(\vvv{K}_k)$. We want to show that the set $g^{-1}(u,i)$ is acyclic. Notice that the set $g^{-1}(u, i)$ is equal to the set $\{(u,f^{-1}(i)) \in V(G) \times V(H)\}$. Since $f$ is an acyclic homomorphism, the subdigraph of $H$ induced by the set $f^{-1}(i)$ is acyclic. Now, note that the subdigraph induced by  $g^{-1}(u, i)$ is a copy of the subdigraph of $H$ induced by $f^{-1}(i)$. It thus follows that the set $g^{-1}(u, i)$ is acyclic, and that $g$ is an acyclic homomorphism, as we wanted to prove.

\item[$(\geq)$]

Now, let $f : V(G[H]) \rightarrow [\ell]$ be and acyclic coloring of $G[H]$. We want to build an acyclic coloring of $G[\vvv{K}_k]$ with at most $\ell$ colors using $f$, thus getting that $\cn(G[\vvv{K}])\le \cn(G[H])$ and finishing our proof. 
For each $u \in V(G)$ denote by $H_u $ the copy of $H$ related to $u$, i.e., the subdigraph of $G[H]$ induced by $\{u\} \times V(H)$. 
Let $C^u = \{c^u_1,..., c^u_{\ell_u} \} \subseteq [\ell]$ be the set of colors used in $H_u$, for each $u\in V(G)$; more formally, $C^u = f(V(H_u))$. Because $f$ restricted to $H_u$ is an acyclic coloring of $H_u$, we get that $\ell_u\ge k$. 
Denoting the vertex set of $V(\vvv{K}_k)$ by $[k]$, consider $g: V(G[\vvv{K}_k]) \rightarrow [\ell]$ such that $g(u, i)=c^u_i$ for each $i\in [k]$ and $u\in V(G)$.  We show that $g$ is an acyclic coloring of $G[\vvv{K}_k]$. 
Suppose, by contradiction, that $C=(w_1, \ldots, w_t)$ is a monochromatic dicycle of color $c$ and write $w_j=(u_j, i_j) \in V(G[\vvv{K}_k])$, for each $j \in [t]$. 
Let $U=\{u_1, \ldots, u_t\}$, and notice that we must have $|U|\ge 2$, since there is at most one vertex of color $c$ in $u_j[\vvv{K}_k]$ for every $j\in [t]$. 
Thus, by Lemma \ref{lemma:key} we get that $U$  contains a dicycle of $G$. Now, for each $u_j$ with $j \in [t]$, we pick $x_{j} \in V(H)$ such that  $f(u_j,x_j)= c$; it exists since $c\in C^u$ and by the definition of the set of colors $C^u$. Denote $(u_j,x_j)$ by $w'_j$, for each $j\in [t]$, and note that, since $U$ contains a dicycle, the  subdigraph of $G[H]$ induced by the set $W=\{w'_{1}, \ldots, w'_{t}\}$ contains a monochromatic dicycle, contradicting the fact that $f$ is an acyclic coloring of $G[H]$. Therefore, $g$ is an acyclic coloring of $G[\vvv{K}_k]$, as we wanted to prove.
\end{itemize}
\end{proof}


\section{Products of Dicycles}\label{sec:dicycles}

In this section, we present the proofs of the results involving dicycles, more especifically, we prove Theorems~\ref{thm:positiveAnswer},~\ref{thm:ProdForte_outros},~\ref{thm:prodforte} and~\ref{thm:nlechi}.

To prove Theorem~\ref{thm:positiveAnswer}, we first prove the following lemma.

\begin{lemma}\label{lem:directCycles}
Let $\vv{C}_n$, $\vv{C}_m$ be dicycles on $n,m$ vertices, respectively. Then, $$\cn(\vv{C}_n \times \vv{C}_m)=2.$$
\end{lemma}
\begin{proof} Write $V(\vv{C}_{n}) = [n]$ and $V(\vv{C}_{m}) = [m]$. By the definition of direct product, each vertex $(i, j) \in V(\vv{C}_{n} \times \vv{C}_{m})$ has in-degree, as well as out-degree, equal to $1$. 
It follows that each strong component of $\vv{C}_{n} \times \vv{C}_{m}$ is a dicycle; hence, we need at least~2 colors in any acyclic coloring of $\vv{C}_{n} \times \vv{C}_{m}$. Equality follows by Theorem~\ref{thm:folkore}.
\end{proof}

Below, we restate Theorem~\ref{thm:positiveAnswer}.

\thmpositiveAnswer*
\begin{proof}
    When $\cn(G) = 1$, we have equality already by Corollary~\ref{cor:DAG}; so suppose $\cn(G) = 2$. Since $\cn(H)\ge \cn(G)$, we get that both $G$ and $H$ contain dicycles, say $C_G\subseteq G$ and $C_H\subseteq H$. By Theorem~\ref{thm:folkore}, we know that $\cn(G\times H)\le 2$, and because $C_G\times C_H\subseteq G\times H$, it follows that $\cn(G\times H)\ge 2$ by the above lemma.
\end{proof}

\thmprodforteoutros*
\begin{proof}
We first prove that $\chi(H\boxtimes C_3) = 6$. To see that $\chi(H \boxtimes C_{3}) \ge 6$, note that $K_6 \cong K_2 \boxtimes C_3 \subseteq H \boxtimes C_{3}$. As for $\chi(H \boxtimes C_{3}) \le 6$, it follows from the the fact that $H\boxtimes C_3\subseteq H[C_3]$ and from Theorem~\ref{prop:produto}.

    In \cite{V.79} the author proves that if the two factors have at least one edge, then the strong product has at least the maximum between the chromatic numbers of the factors plus $2$. In this case, $\chi(H \boxtimes C_{2n+1}) \ge 3 + 2 = 5$. To prove that $\chi(H\boxtimes C_{2n+1}) \le 5$, we present a $5$-coloring of $H \boxtimes C_{2n+1}$. Let ${A,B}$ be the bipartition of $H$. Regarding the copies of $C_{2n+1}$ with respect to the vertices in $A$, we color the first three vertices with colors $1,2,3$ respectively and alternate colors $4,5$ in each one of the remaining vertices, and regarding the copies of $C_{2n+1}$ with respect to the vertices in $B$, we color the first three vertices with colors $3,4,5$ respectively and alternate colors $1,2$ in each one of the remaining vertices. More formally, write $C_{2n+1}$ as $(x_1,\ldots, x_{2n+1})$ and define $h: H \boxtimes C_{2n+1} \rightarrow [5]$ as follows.
    
    \hspace{-8.4pt}$h(u,x_j)=\begin{cases} 1, \textrm{if $u \in A$ and $j=1$, or $u \in B$ and $j=2q$, $q \in \{2,\ldots,n\}$;}\\
    2, \textrm{if $u \in A$ and $j=2$, or $u \in B$ and $j=2q+1$, $q \in \{2,\ldots,n\}$;}\\
    3, \textrm{if $u \in A$ and $j=3$, or $u \in B$ and $j=1$;}\\
    4, \textrm{if $u \in A$ and $j=2q$, $q \in \{2,\ldots,n\}$, or $u \in B$ and $j=2$;}\\
    5, \textrm{if $u \in A$ and $j=2q+1$, $q \in \{2,\ldots,n\}$, or $u \in B$ and $j=3$.}

    \end{cases}$

    We want to show that every color class is an independent set. Recall that, by the definition of strong product (taking the indices mod $2n+1$), we have
    \begin{equation*}
        N(u,x_j) = \left(N(u) \times \{x_{j-1},x_j,x_{j+1}\}\right) \cup \left( \{u\}\times \{x_{j-1},x_{j+1} \}\right))
    \end{equation*}
    
    Observe that $h(V(H)\times \{x_j\})$ has exactly two colors, and that these colors are distinct from the colors used in $V(H)\times \{x_{j-1}\}$ and from the colors used in $V(H)\times \{x_{j+1}\}$. Because these colors also respect the bipartition of $H$, we get the desired proper coloring of $H\boxtimes C_{2n+1}$.
    
\end{proof}

Now, we investigate the dichromatic number of the strong product of two dicycles. 

\thmprodforte*
\begin{proof}
For Item~1, just observe that $\vv{C}_2 \boxtimes \vv{C}_2 = \vvv{K}_4$; hence $\cn(\vv{C}_2 \boxtimes \vv{C}_2) = \cn(\vvv{K}_4)=4$. Now, we prove Item~2. First consider $n=2$. Recall that $\vv{C}_3 \boxtimes \vv{C}_2 \subseteq \vv{C}_3[\vv{C}_2]$,  and since we will see $\cn(\vv{C}_3[\vv{C}_2]) = 3$ in Theorem \ref{thm:nlechi}, it follows that $\cn(\vv{C}_3 \boxtimes \vv{C}_2) \leq 3$.  To prove that $\cn(\vv{C}_3 \boxtimes \vv{C}_2) \ge 3$, consider a 2-coloring $f:V(\vv{C}_3) \times V(\vv{C}_2) \rightarrow [2]$. To follow the demonstration, see Figure \ref{fig:3x2}. Write  $V(\vv{C}_3 \boxtimes \vv{C}_2)=\{(u_i, x_j)|i \in [3],j\in [2]\}$. Since $((u_1,x_1), (u_2,x_1), (u_3,x_1))$ is a dicycle, we must have two colors appearing in these vertices. Furthermore, by the Pigeon Hole Principle, some color appears in exactly two vertices. Suppose, without loss of generality, that $f(u_1,x_1)=f(u_2,x_1)=1$. Now, for every $i \in [3]$, because $((u_i,x_1), (u_i,x_2))$ is a dicycle, we have $f(u_i,x_1)\neq f(u_i,x_2)$. Observe that this gives us exactly the coloring presented in Figure~\ref{fig:3x2}, a contradiction since $((u_1,x_1), (u_2,x_1), (u_3,x_2))$ is a monochromatic dicycle of color $1$. 

\iffigures{
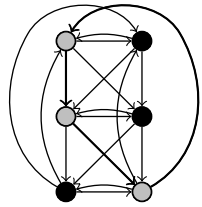
\begin{figure}[H]
    \begin{center}
\captionsetup{justification=centering}
         \begin{tikzpicture}[scale=1]
 \pgfsetlinewidth{.5pt}

  \tikzstyle{blackvertex}=[draw,circle,fill=black,minimum size=7pt,inner sep=0pt]
   \tikzstyle{vertex}=[draw,thick,circle,fill=white,minimum size=7pt,inner sep=0pt]
   \tikzstyle{grayvertex}=[draw,circle,fill=black!25,minimum size=7pt,inner sep=0pt]

    \foreach \i in {0,1}{
    \foreach \j in {0,1,2}{
  \node[vertex] (\i\j) at (\i,\j){};
    }
  }
    \foreach \i in {0,1}{
 
  \draw[->] (\i1) -- (\i0);
  \draw[->] (\i2) -- (\i1);
  \draw[->] (\i0)to[out=120,in=-120] (\i2);
  }
  
  \foreach \j in {0,1,2}{
  \draw[->] (0\j) -- (1\j);

  \draw[->] (1\j)to[out=165,in=15] (0\j);
  
  }
  \draw[->] (00)to[out=150,in=240] (-0.5,2.1) to[out=60,in=120] (12);
  \draw[->] (10)to[out=30,in=-60] (1.5,2.1) to[out=120,in=60] (02);
  
  \draw[->] (02) -- (11);
  \draw[->] (12) -- (01);
  \draw[->] (01) -- (10);
  \draw[->] (11) -- (00);
  
    \node[grayvertex] (02) at (0,2){};
  \node[grayvertex] (01) at (0,1){};
  \node[grayvertex] (10) at (1,0){};
  \node[blackvertex] (00) at (0,0){};
  \node[blackvertex] (12) at (1,2){};
  \node[blackvertex] (11) at (1,1){};
  
  \draw[black,thick,->] (02) -- (01);
  \draw[black,thick,->] (01) -- (10);
  \draw[black,thick,->] (10)to[out=30,in=-60] (1.5,2.1) to[out=120,in=60] (02);
     \end{tikzpicture}
            \caption{2-coloring of $C_3\boxtimes C_2$ constructed in the proof of Theorem~\ref{thm:prodforte}. Color~1 is represented in gray, and color~2 in black.}
             \label{fig:3x2}
    \end{center}
\end{figure}}\fi
\smallskip

Now consider $n=3$. Again, we have $ \cn(\vv{C}_3 \boxtimes \vv{C}_3) \leq 3 $, because $\vv{C}_3 \boxtimes \vv{C}_3 \subseteq \vv{C}_3[\vv{C}_3] $ and by the Theorem~\ref{thm:nlechi}. Now, we show that $\cn(\vv{C}_3 \boxtimes \vv{C}_3) > 2$ by testing each possible $2$-coloring of $\vv{C}_3 \boxtimes \vv{C}_3$. Note that, up to relabeling, there exists only one possible 2-coloring of the first column, $C^1$. We use the colors gray $(G)$ and black  $(B)$ in the colorings in Figure \ref{fig:C3xC3_colorings}. Without loss of generality, we fix the coloring of $C^1$ into $GGB$. Then, again because the first row $L^1$ cannot be monochromatic, it can only be colored in two possible ways: $GBB$ or $GGB$ (note that $GBG$ is also a possibility, but it is only a cyclic permutation of $GGB$, in other words, it is also a relabeling). Observe that, for each choice of coloring of $L^1$, there are 16 possible colorings of the remaining vertices. We exhibit the $16$ possible colorings for each coloring of $L^1$, and in each of them, we present a monochromatic dicycle (emphasized arcs in the figure). Therefore, $\cn(\vv{C}_3 \boxtimes \vv{C}_3) > 2$.

\iffigures{
\begin{figure}[H]
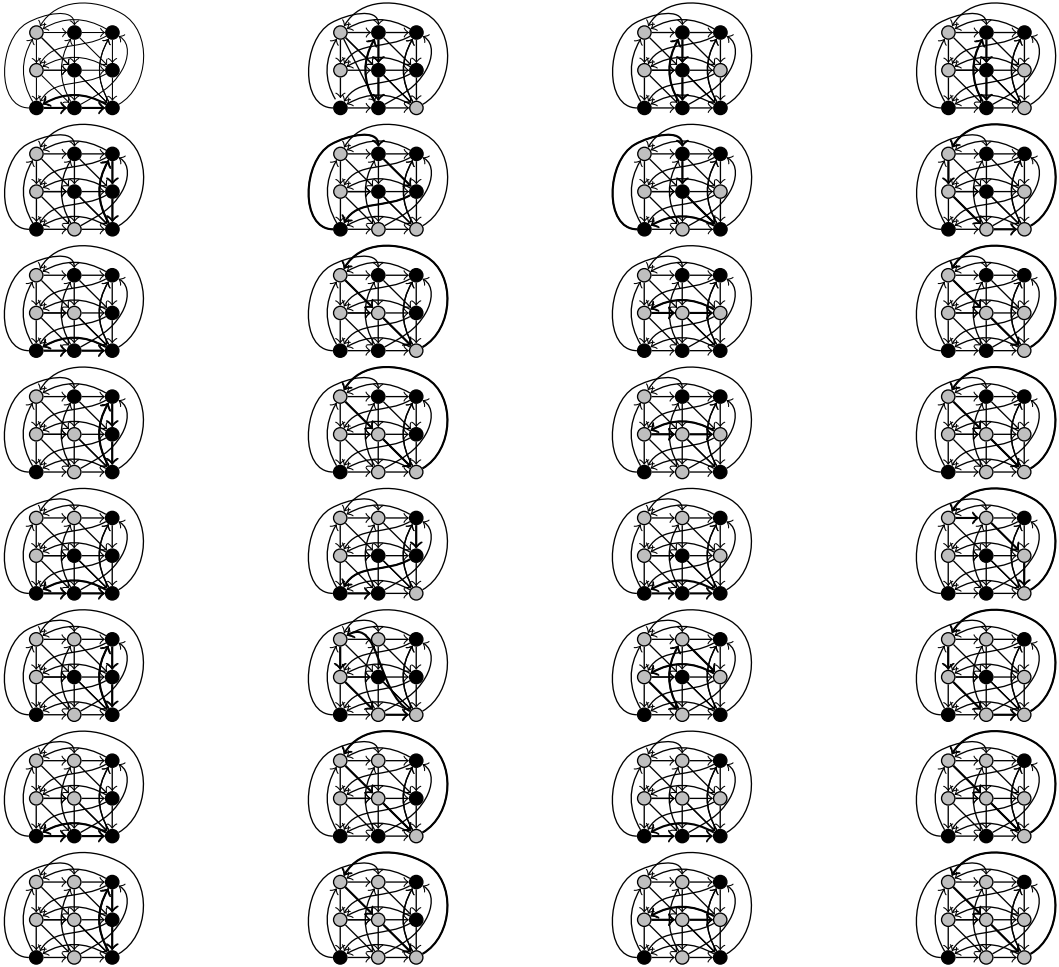

    \begin{center}
        \captionsetup{justification=centering}

 \hspace{-12cm}


                \caption{The $16$ colorings of $\vv{C}_3 \boxtimes \vv{C}_3$, fixing the coloring of $C^1$ in $GGB$  and the coloring of $L^1$ in $GBB$ , and the $16$ colorings of $\vv{C}_3 \boxtimes \vv{C}_3$, fixing the coloring of $C^1$ in $GGB$  and the coloring of $L^1$ in $GGB$. In each of them we highlight the edges forming a dicycle.}
                \label{fig:C3xC3_colorings}
    \end{center}
\end{figure}

\smallskip
}\fi

Finally, we prove item 3. Since $\vv{C}_m\times \vv{C}_n\subseteq \vv{C}_m\boxtimes \vv{C}_n$, by Lemma~\ref{lem:directCycles} we get $ \cn(\vv{C}_m \boxtimes \vv{C}_n) \ge 2$. Thus, in what follows we simply provide acyclic colorings with~2 colors. Write $\vv{C}_m = (u_1,\cdots, u_m)$ and $\vv{C}_n = (x_1,\cdots, x_n)$. For each $i\in [m]$, we define the \emph{$i$-th row} of $\vv{C}_m\boxtimes\vv{C}_n$ as the subset $R^i = \{(u_i,x_j)\mid j\in [n]\}$; similarly, for each $j\in [n]$, the \emph{$j$-th column} of $\vv{C}_m \boxtimes\vv{C}_n$ is the subset $C^j = \{(u_i,x_j)\mid i\in [m]\}$. 
Define the following coloring of $\vv{C}_m\boxtimes \vv{C}_n$ (observe Figure~\ref{fig:mxn}):
\[
f(u_i,x_j) = \left\{\begin{array}{cl}
     1 &  \text{, if either $(i,j)\in \{(2,3),(3,2),(3,3)\}$,}\\&\text{ or $i=1$ and $j\ge 3$, or $j=1$ and $i\ge 3$;} \\
     2 & \text{, otherwise}
\end{array}\right.
\]

\iffigures{
\begin{figure}[H]
    \begin{center}
    \captionsetup{justification=centering}

 \begin{tikzpicture}[scale=0.6]
 
  \pgfsetlinewidth{.7pt}

  \tikzstyle{blackvertex}=[draw,circle,fill=black,minimum size=7pt,inner sep=0pt]
   \tikzstyle{vertex}=[draw,thick,circle,fill=white,minimum size=7pt,inner sep=0pt]
   \tikzstyle{grayvertex}=[draw,circle,fill=black!25,minimum size=7pt,inner sep=0pt] 
  
    \foreach \i in {0,1,2,3,4,5}{
        \foreach \j in {0,1,2,3,4,5}{
            \node[vertex] (\i\j) at (\i,\j){};
        }
    }
    
    \foreach \i in {0,1,2,3,4,5}{
        \draw[->] (\i5) -- (\i4);
        \draw[->] (\i4) -- (\i3);
        \draw[->] (\i3) -- (\i2);
        \draw[->] (\i2) -- (\i1);
        \draw[dotted] (\i1) -- (\i0);
        \draw[->] (\i0)to[out=120,in=-120] (\i5);
    }

  \foreach \j in {0,1,2,3,4,5}{
    \draw[->] (0\j) -- (1\j);
    \draw[->] (1\j) -- (2\j);
    \draw[->] (2\j) -- (3\j);
    \draw[->] (3\j) -- (4\j);
    \draw[dotted] (4\j) -- (5\j);
    \draw[->] (5\j)to[out=150,in=30] (0\j);
  }
  \node[grayvertex] (0x) at (0,5){};
  \node[grayvertex] (Ax) at (0,4){};
  \node[blackvertex] (Bx) at (0,3){};
  \node[blackvertex] (Cx) at (0,2){};
  \node[blackvertex] (Dx) at (0,1){};
  \node[blackvertex] (Ex) at (0,0){};
  \node[grayvertex] (0y) at (1,5){};
  \node[grayvertex] (Ay) at (1,4){};
  \node[blackvertex] (By) at (1,3){};
  \node[grayvertex] (Cy) at (1,2){};
  \node[grayvertex] (Dy) at (1,1){};
  \node[grayvertex] (Ey) at (1,0){};
  \node[blackvertex] (0z) at (2,5){};
  \node[blackvertex] (Az) at (2,4){};
  \node[blackvertex] (Bz) at (2,3){};
  \node[grayvertex] (Cz) at (2,2){};
  \node[grayvertex] (Dz) at (2,1){};
  \node[grayvertex] (Ez) at (2,0){};
  \node[blackvertex] (0w) at (3,5){};
  \node[grayvertex] (Aw) at (3,4){};
  \node[grayvertex] (Bw) at (3,3){};
  \node[grayvertex] (Cw) at (3,2){};
  \node[grayvertex] (Dw) at (3,1){};
  \node[grayvertex] (Ew) at (3,0){};
  \node[blackvertex] (0R) at (4,5) {};

  \node[grayvertex] (AR) at (4,0) {};
  \node[grayvertex] (BR) at (4,1) {};
  \node[grayvertex] (CR) at (4,2) {};
  \node[grayvertex] (DR) at (4,3) {};
  \node[grayvertex] (ER) at (4,4) {};
  \node[blackvertex] (0S) at (5,5) {};
  \node[grayvertex] (AS) at (5,0) {};
  \node[grayvertex] (BS) at (5,1) {};
  \node[grayvertex] (CS) at (5,2) {};
  \node[grayvertex] (DS) at (5,3) {};
  \node[grayvertex] (ES) at (5,4) {};
 
 \end{tikzpicture}

\caption{Acyclic coloring using~2 colors of $\vv{C}_m\boxtimes \vv{C}_n$, with $m\ge 4$ and $n \ge 4$. Here we omit diagonal edges to keep the figure clearer. In each row, we have a copy of $\vv{C}_n$, and in each column, a copy of $\vv{C}_m$. Also, color 1 is presented in black, while color~2, in gray.}
             \label{fig:mxn}
    \end{center}
\end{figure}
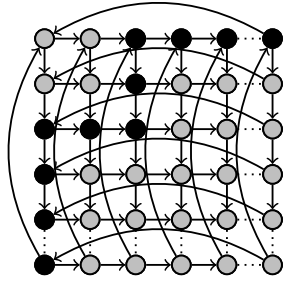}\fi

First, we prove that color~2 cannot contain a dicycle. Suppose otherwise and let $C = (w_1,\cdots, w_t)$ be such a cycle. Also, write $w_\ell = (u_{i_\ell},x_{j_\ell})$ for each $\ell\in [t]$, and let $U = \{u_{i_1},\cdots, u_{i_t}\}$ and $X = \{x_{j_1},\cdots,x_{j_t}\}$. 
Observe that Lemma~\ref{lemma:cycleinColumn} tells us that either $U = V(\vv{C}_m)$ or $X = V(\vv{C}_n)$. If the former occurs (i.e., $C$ contains a vertex of color~2 of each row $R^i$), then, because $(u_1,x_1)$ and $(u_1,x_2)$ are the only vertices in the 1st row of color~2, we get that at least one of them is contained in $C$. Similarly, if the latter occurs, because $(u_1,x_1)$ and $(u_2,x_1)$ are the only vertices in the 1st column of color~2, we get that at least one of them is contained in $C$. Therefore, we always get that $V(C)\cap T\neq \emptyset$, where $T = \{(u_1,x_1),(u_1,x_2),(u_2,x_1)\}$. Now, let $S = T\cup \{(u_2,x_2)\}$, and suppose, without loss of generality, that $w_1\in V(C)\cap T$. Because $N^+(T)\cap f^{-1}(2) \subseteq S$, $N^+(S)\subseteq f^{-1}(1)$ and $S$ is acyclic, we get a contradiction.

Now, we prove that the subdigraph induced by color~1 is acyclic. For this, we start by letting $F$ be the set of vertices colored with~1, and iteratively remove vertices from $F$ when we can conclude that they cannot be contained in any cycle. If $F$ ever becomes acyclic, we are done. Denote by $G$ the subdigraph of $\vv{C}_m\boxtimes\vv{C}_n$ induced by $F$. Observe that if $w\in F$ is either a source or a sink in $G$, then $w$ cannot be contained in any cycle of $G$. Now, observe that the definition of strong product tells us that:
\smallskip
\begin{enumerate}
    \item[(*)] $(u_i,x_j)(u_{i'},x_{j'})\in E(\vv{C}_m\boxtimes \vv{C}_n)$ if, and only if, $i=i'$ and $j' = (j+1)\mod n$, or $i'=(i+1)\mod m$ and $j=j'$, or $i'=(i+1)\mod m$ and $j' = (j+1)\mod n$.
\end{enumerate} 
\medskip

Note that this gives us that $N^-(u_3,x_1) = \{(u_2,x_1),(u_3,x_n),(u_2,x_n)\}$; because all of these vertices are of color~2, we get that $(u_3,x_1)$ is a source in $G$ and therefore can be removed from $F$. Applying similar arguments, we can remove from $F$ the following vertices, in the presented order:

\medskip
\begin{enumerate}
    \item $(u_1,x_3)$: this is because $N^-(u_1,x_3) = \{(u_1,x_2),(u_m,x_3),(u_m,x_2)\}\subseteq f^{-1}(2)$; 
    \item $(u_2,x_3)$: this is because $N^-(u_2,x_3) = \{(u_2,x_2), (u_1,x_3),(u_1,x_2)\}$, and the fact that $(u_1,x_3)$, the only neighbor of $(u_2,x_3)$ of color~1, was already removed from $F$;
    \item $(u_3,x_2)$: this is because $N^-(u_3,x_2) = \{(u_3,x_1),(u_2,x_2),(u_2,x_1)\}$, and the fact that $(u_3,x_1)$ was already removed from $F$;
    \item $(u_3,x_3)$: this is because $N^-(u_3,x_3) = \{(u_3,x_2),(u_2,x_3),(u_2,x_2)\}$, and the fact that $(u_3,x_2)$ and $(u_2,x_3)$ were already removed from $F$.
\end{enumerate}
\medskip

Finally, we get that the remaining vertices in $F$ are $\{(u_i,x_j)\mid \text{ either }i = 1\text{ and }j\in\{4,\cdots,n\}\text{, or }i\in\{4,\cdots,m\}\text{ and }j=1\}$. We prove that the following ordering is a topological ordering of $F$, finishing the proof since we then get that $F$ is acyclic: \[((u_4,x_1),(u_5,x_1),\cdots,(u_m,x_1),(u_1,x_4),(u_1,x_5),\cdots,(u_1,x_n)).\]
Suppose otherwise and let $(u_i,x_j)(u_{i'},x_{j'})\in E(C_m\boxtimes C_n)$ be an edge that does not respect the above ordering (i.e., it goes from a bigger to a smaller vertex in the ordering). 
If $i\neq i'$, observe that either $i=1$ and $i'\in \{4,\ldots,m\}$, or $4\le i' < i \le m$. 
In any case, we get that $i' \neq  (i+1)\mod m$, contradicting $(*)$. Therefore, we must have $i=i'=1$. But now we get $j' < j\le n$, which in turn implies that $j' \neq (j+1)\mod n$, again contradicting $(*)$.
\end{proof}

\medskip

As discussed in Section~\ref{sec:prelim}, the following is a generalization of a result presented in~\cite{PP.16}.

\thmnlechi*

\begin{proof}
In what follows, given a function $f :V(G) \to \mathbb{N}$ and a subdigraph $G' \subseteq G$, we make an abuse of notation and we write $f(G')$ to denote $f(V(G'))$.

Write $V(\vv{C}_n)=\{u_1,\ldots,u_n\}$ and, for each $i \in [n]$, denote by $H_i$ the copy of $H$ related to the vertex $u_i$, in other words, $H_i$ is the subdigraph of $\vv{C}_n[H]$ induced by the set $\{(u_i,x) \mid x \in V(H)\}$. We first show an useful fact. Consider a coloring $f:V(\vv{C}_n[H])\to[\ell]$ such that $f$ restricted to $H_i$ is acyclic, for each $i \in [n]$. Below we show that:
\begin{itemize}
         \item [$(*)$]\label{ast} $f$ is an acyclic coloring of $\vv{C}_n[H]$ if, and only if, for each $c \in [\ell]$, there exists $i \in [n]$ such that $c \notin f(H_i)$.
        \end{itemize}
        
        In other words, no color  $c \in [\ell]$ is such that $c$ appears in $H_i$ for every $i \in [n]$. 
        We first prove that the condition is necessary. For this, let $c \in [\ell]$ be such that $c \in f(H_i)$ for every $i \in [n]$. Then, for each $i \in [n]$, pick a vertex $x_i \in V(H)$ such that $f(u_i,x_i)=c$. By the definition of  lexicographic product, we get that $((u_1,x_1),\ldots,(u_n,x_n))$ is a monochromatic dicycle of color $c$; hence $f$ is not an acyclic coloring of $\vv{C}_n[H]$. 
        Now, to prove sufficiency, suppose that $f$ yields a monochromatic dicycle $C$ of color $c \in [\ell]$ in $\vv{C}_n[H]$. Write $C=(w_1,\ldots,w_t)$ and  $w_j = (u_{i_j},x_j) \in V(\vv{C}_n[H])$, for each $j \in [t]$. Also, let $U = \{u_{i_1},\cdots,u_{i_t}\}$. Since $f$ restricted to $H_i$ is acyclic, for each $i \in [n]$, we cannot have $|U|=1$. By Lemma \ref{lemma:key}, it follows that $U$ contains a dicycle. Since $\vv{C}_n$ is itself a dicycle, we can only have $U = V(\vv{C}_n)$. That way, for every $i \in [n]$, there exists a vertex $x_j \in V(H)$ such that $f(u_i,x_j)=c$, and $c$ is the desired color.
        
        Now, denote $\cn(H)$ by  $k$, and the value $k+ \lceil \frac{k}{n-1}\rceil$ by $k'$. To show that $\cn(\vv{C}_n[H])=k'$, denote by $s$ the value $k'-k$. First we provide an acyclic coloring for $\vv{C}_n[H]$ with $k'$ colors. To do this, for each $i\in [n]$ we pick the colors that will not appear in $H_i$, then we apply (*). So choose $n$ subsets of $[k']$, $M_1,\ldots, M_n$ , each of cardinality $s$, in a way that $\bigcup_{i=1}^n M_i=[k']$. Notice that, by the Pigeon Hole Principle, for this to be possible it is enough to ensure that $n \ge \frac {k'}{s}$. Indeed, this holds since:
        
        \[\frac{k'}{s}= \frac{k+s}{s}=\frac{k}{s} +1 = \frac{k}{\lceil \frac{k}{n-1}\rceil} + 1 \le \frac{k}{\left(\frac{k}{n-1}\right)} + 1 = n .\]
        
        Now, for each $i \in [n]$, let $f_i$ be a coloring of $H_i$ that uses the colors $[k']\setminus M_i$, and let $f$ be the coloring of $\vv{C}_n[H]$ obtained by the union of these colorings. Since for every $c \in [k']$, there exists $i \in [n]$ such that $c \in M_i$ (which implies that $c \notin f(H_i)$), we have that $f$ satisfies $(*)$ and, therefore, is an acyclic coloring of $\vv{C}_n[H]$.
        
        Finally, let $f$ be an acyclic coloring of $\vv{C}_n[H]$. It remains to show that if $f$ uses $k+ \ell$ colors, then $\ell \ge s$ (i.e., $k+ \ell \ge k'$). For each $i\in [n]$, denote by $f_i$ the coloring $f$ restricted to $H_i$, and  notice that, since $f_i$ must be an acyclic coloring of $H_i$, we have that the set of colors not used in $H_i$ has cardinality at most $\ell$. Denote this set by $M_i$, i.e., let $M_i = [k+\ell] \setminus f(H_i)$ for each $i\in [n]$. Now, by $(*)$ we get that $\bigcup_{i=1}^n M_i = [k+ \ell]$. Therefore, $k+\ell = \left|\bigcup_{i=1}^n M_i\right| \le \ell n$. It follows that $\ell \ge \frac{k}{n-1}$, and since $\ell$ is an integer, we have $\ell \ge s = \lceil \frac{k}{n-1}\rceil$, as we wanted to prove.

\end{proof}

\section{Parameterizing by the Treewidth}\label{sec:treewidth}

Given an undirected graph $G$, a \emph{tree decomposition of $G$} is a pair $(T,{\cal X})$ where $T$ is a tree, and ${\cal X}$ is a family that assigns to each $t\in V(T)$ a subset of $V(G)$, denoted by $X_t$, that satisfies:
\begin{enumerate}
    \item $\bigcup_{t\in V(T)}X_t = V(G)$;
    \item For every $uv\in E(G)$, there exists $t\in V(T)$ such that $\{u,v\}\subseteq X_t$; and
    \item For every $t,t'\in V(T)$ and every $t''$ in the $t,t'$-path in $T$, we have that $X_t\cap X_{t'}\subseteq X_{t''}$.
\end{enumerate}

The \emph{width of $(T,{\cal X})$} is given by $\max_{t\in V(T)}|X_t| - 1$, while the \emph{treewidth of $G$} is the minimum width over all tree decompositions of $G$; it is denoted by $tw(G)$. Also, the treewidth of a digraph is equal to the treewidth of its underlying graph. We mention that there are definitions of directed tree decomposition, but that these usually do not lead to $\FPT$ algorithms~\cite{JRST.01,R.99}.
The main result of this section is an algorithm {\FPT}  when parameterized by the treewidth that computes the dichromatic number of a given digraph. 
The {\FPT} time complexity is obtained thanks to the largely known fact that $\chi(H)\le tw(H)+1$ for every graph $H$ (folklore), and because $\cn(G)\le \chi(U(G))$ for every digraph $G$. Before we present our algorithm, we show a slightly better upper bound for oriented graphs. 

\begin{theorem}
\label{thm:tworient} Let $G$ be a graph and $D$ be an orientation of $G$. Then, \[\cn(D) \le \left\lceil \frac{tw(G)+1}{2} \right \rceil.\]
\end{theorem}
\begin{proof}
By induction on the number of vertices. Let $(T,\mathcal{X})$ be a nice tree decomposition of $G$ with width $tw(G)$ and $\ell$ be a leaf node of $T$. We can suppose, without loss of generality, that there exists $u\in X_\ell$ such that $X_\ell$ is the only subset containing $u$. Indeed, if this is not the case, then $X_\ell \subseteq X_t$, where $t$ is the neighbor of $\ell$ in $T$, and we can simply remove $\ell$ from the tree decomposition. Now, let $G'=G-u$ and $D'=D-u$. By induction hypothesis, there exists an acyclic coloring of $D'$ with at most $k' = \left\lceil \frac{tw(G')+1}{2} \right \rceil$ colors. Observe that if there exists a color $c$ such that $u$ has at most one neighbor colored with $c$, then we can color $u$ with $c$ without creating monochromatic dicycles. If $tw(G')$ is even, then $k'> \frac{tw(G')+1}{2}$, i.e., $|X_\ell\setminus \{u\}|\le tw(G')+1 < 2k'$, which by the Pigeonhole Principle and the fact that $N(u)\subset X_\ell$ implies that there exists the desired color. Otherwise, write $tw(G')$ as $2t+1$ and suppose that $u$ has at least two neighbors in each color, as otherwise again we can extend the coloring.  Because $N(u)\subseteq X_\ell$, we get $|X_\ell| \ge 2k'+1 \ge tw(G')+2$. Since $|X_\ell| \le tw(G)+1\le tw(G')+2$, we get equality and $k = \left\lceil \frac{tw(G)+1}{2} \right \rceil = \left\lceil \frac{2t+3}{2} \right \rceil = t+2 = k'+1$. Hence we can use a new color on $u$. 
\end{proof}

We now present our {\FPT} algorithm. As usual, it is a dynamic programming algorithm that uses a particular type of tree decomposition, that we define next.

In order to avoid confusion, the vertices of $T$ are referred to as \emph{nodes}. A \emph{nice tree decomposition} is a tree decomposition $(T,{\cal X})$, where $T$ is rooted in a vertex $r$, and each $t\in V(T)$ is one of the following types of nodes:
\begin{itemize}
    \item \emph{Leaf node}: in this case, $t$ is a leaf of $T$ and $X_t$ is unitary;
    \item \emph{Forget node}: $t$ has exactly one child, $t'$, and $X_t = X_{t'}\setminus \{v\}$. We say that $v$ is \emph{forgotten by $t$}; 
    \item \emph{Introduce node}: $t$ has exactly one child, $t'$, and $X_t = X_{t'}\cup \{v\}$. We say that $v$ is \emph{introduced by $t$}; and
    \item \emph{Join node}: $t$ has exactly two children, $t_1,t_2$, and $X_t = X_{t_1} = X_{t_2}$.
\end{itemize}

It is largely known that, in time $2^\mathcal{O}(k)\cdot n$, one can either conclude that the treewidth of $G$ is larger than $k$, or obtain a tree decomposition of $G$ of width at most $5k + 4$~\cite{B.16}. Also, given a tree decomposition $\mathcal{T}=(T,\mathcal{X})$ of $G$ of width at most $k$, one can in time $\mathcal{O}(k^2 \cdot \max (|V(T)|,|V(G)|))$ compute a nice tree decomposition of $G$ of width at most $k$ that has at most $\mathcal{O}(k|V(G)|)$ nodes \cite{CMF.15}. We consider that we know such a tree decomposition and later discuss about the final complexity.

Given a node $t$ of the considered tree decomposition, we define $G_t$ as the subgraph of $G$ induced by the vertices in $X_t$ along with the vertices in all the subsets $X_t'$ such that $t'$ is a descendant of $t$ in $T$.

\begin{figure}[H]
\begin{center}

\begin{tikzpicture}
 
  \pgfsetlinewidth{.7pt}

  \tikzstyle{blackvertex}=[draw,circle,fill=black,minimum size=7pt,inner sep=0pt]
   \tikzstyle{vertex}=[draw,thick,circle,fill=white,minimum size=7pt,inner sep=0pt]
   \tikzstyle{grayvertex}=[draw,circle,fill=black!25,minimum size=10pt,inner sep=0pt]

 \node at (0,3.5) (z){};
 \node[draw,circle,label=left:{\{$v_2,v_3,v_4$\}},label=right:$t$] at (0,2.7) (t){};
 \draw[dotted] (t) to (z);
 \node[draw,circle,label=left:{$\{v_2,v_4\}$},label=right:$t'$] at (0,1.5) (t'){};
 \draw (t) -- (t');
 \node[draw,circle,label=left:{$\{v_2,v_4,v_5\}$}] at (0,0.3) (w){};
 \draw (t') -- (w);
 \node at (0,-.5) (a){};
 \draw [dotted] (w) to (a);

 \node[grayvertex] (11) at (3,0) {$v_1$};
 \node[grayvertex] (22) at (3,2) {$v_2$};
 \node[grayvertex] (33) at (4,3) {$v_3$};
 \node[grayvertex] (44) at (5,2) {$v_4$};
 \node[vertex] (55) at (5,0) {$v_5$};
 
 \draw[->] (11)--(22);
 \draw[->] (22)--(33);
 \draw[->] (33)--(44);
 \draw[->] (44)--(55);
 \draw[->] (55)--(11);
 \draw[->] (44)--(11);

 \draw[->] (55)--(22);
 \draw[->] (11) to[out=-30,in=-150] (55);

 \end{tikzpicture}

\end{center}
\caption{To the left, we represent part of a nice tree decomposition of the digraph depicted to the right. The coloring of $X_t$, an introduce node, cannot be extended to an acyclic 2-coloring of $G_t$, even though there exists an acyclic coloring that extends the 2-coloring of $X_{t'}$ and the coloring constrained to $X_t$ is an acyclic coloring.}
\label{fig:coincidir} 
\end{figure}
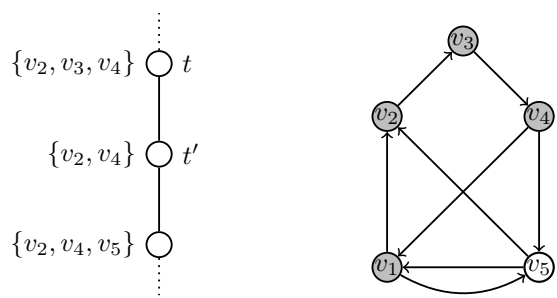

Now, given a digraph $G$, we work on a nice tree decomposition $(T,{\cal X})$ of the underlying graph of $G$, computing information for the children of each node $t$ before being able to compute information for $t$. In order to do that, in the classical proper coloring problem, it is enough to keep, for every possible coloring $f$ of $X_t$, whether $f$ can be extended to $G_t$. This is because $X_t$ is a separator of $G$, for each node $t$; hence, when introducing a new vertex $v$ or joining two colorings on a join node, it is enough to know that $f$ can be extended to the rest of the graph. However, for acyclic colorings, it can happen that, when adding a new vertex or joining two colorings on a join node, a cycle $C$ is created within a coloring, with $V(C)\nsubseteq X_t$. See for instance the coloring of Figure \ref{fig:coincidir}.  This tells us that it is not enough to know whether an acyclic coloring $f$ can be extended to an acyclic coloring $f'$ of $G_t$, we also need to know whether there are monochromatic paths between pairs of nodes of $X_t$ in $f'$. For this, we define the notion of representation.

Given a node $t\in V(T)$ and a $k$-coloring $f'$ of $G_t$, we say that $f'$ is  \emph{represented} by $H=(X_t, E')$ in $X_t$ if for every $u,v \in X_t$, $u\neq v$, we have that $uv \in E'$ if, and only if, there is a monochromatic $(u,v)$-path by $f'$ in $G_t$. Let also $f$ be a $k$-coloring of $X_t$. A digraph $H$ is said to be \emph{transitive} if whenever $uv$ and $vw$ are arcs of $H$, we get that $uw$ is also an arc of $H$. 
The \emph{transitive closure} of a digraph $H$ is the minimal transitive digraph containing $H$. Let $f$ be a coloring of $X_t$, and $H' = (X_t,E')$ contain all the monochromatic arcs of $G$, i.e., $uv\in E(H')$ if and only if $uv\in E(G[X_t])$ and $f(u) = f(v)$. 
The \emph{minimal representation of $f$} is the transitive closure of $H'$; we denote it by $H_f$. 
The following simple proposition will allow us to compute trivial entries of our table.

\begin{proposition}\label{prop:trivialzero}
Let $G$ be a digraph, $(T,{\cal X})$ be a tree decomposition of the underlying graph of $G$, and $t$ be a node of $T$. If $f'$ is an acyclic coloring of $G_t$ represented by $H$ in $X_t$, then $H$ is acyclic, transitive and $H_f\subseteq H$, where $f$ is equal to $f'$ restricted to $X_t$.
\end{proposition}
\begin{proof}
Since each arc $uv$ in $H$ is related to a monochromatic $(u,v)$-path in $G_t$, we get that $H$ is acyclic as otherwise we would get a monochromatic cycle in $f'$. Also, if $(u,v),(v,w)\in E(H)$, then there is a monochromatic  $(u,v)$-path and a monochromatic $(v,w)$-path in $f'$, which implies in the existence of a monochromatic $(u,w)$-path in $f'$, and hence $(u,w)\in E(H)$. Finally, since an arc is itself a path, we get that $H$ must contain the monochromatic arcs of $G$, and because $H$ is transitive, it must also contain $H_f$.
\end{proof}

Our table will be indexed by the possible colorings of $X_t$ together with the possible representations of acyclic colorings. So, define the set of $k$-colorings of $X_t$, $\mathcal{F}_t = \{f: X_t \rightarrow [k] \}$, and for each $k$-coloring $f \in \mathcal{F}_t$, define the family $\mathcal{H}_t(f)$  containing all possible digraph with vertex set $X_t$.

Now, for every $f \in \mathcal{F}_t$ and for every $H \in \mathcal{H}_t(f)$, we define:

\[	c_t[f,H]=\begin{cases} 1,& \textrm{if there is an acyclic $k$-coloring $f'$ of $G_t$}\\
& \textrm{which extends $f$ and is represented by $H$ in $X_t$};\\
							0,& \textrm{otherwise.}
			\end{cases}\]

Observe that $G$ has an acyclic $k$-coloring if and only if there exists a $k$-coloring $f$ of $X_r$ and a graph $H\in {\cal H}_r(f)$ such that $c_r[f,H]=1$. 
We now show how to compute $c_t[f,H]$ for each type of node $t$, given that the complete tables of its children are known. We dedicate a lemma to each type of node. First, observe that Proposition~\ref{prop:trivialzero} already tells us that $c_t[f,H]=0$ whenever $H$ does not satisfy the proposition. The following computes the simplest tables, the ones related to leaves of $T$. Because in this case $X_t$ is a single vertex, we get that the only graph in ${\cal H}_t(f)$ is the graph $(X_t,\emptyset)$.

\begin{lemma}\label{lemma:folha}
Let $G$ be a digraph, $(T,{\cal X})$ a nice tree decomposition of $G$, $t$ be a leaf node of $T$, and $ f \in \mathcal{F}_t$. Then ${\cal H}_t(f) = \{H\}$, where $H = (X_t,\emptyset)$, and $c_t[f,H]=1$. Thus, it takes time $O(1)$ to decide whether $c_t[f,H]=1$.
\end{lemma}

\vspace{.5cm}
Now consider $t$ to be a forget node with child $t'$, and let $v\in V(G)$ for which $X_t = X_{t'}\setminus \{v\}$. Consider also $f\in {\cal F}_t$ and $H\in {\cal H}_t(f)$. In order to compute $c_t[f,H]$, we need to analyse the entries of $t'$ containing colorings that extend $f$ to a coloring of $X_{t'}$, i.e, they extend $f$ in order to color $v$. For each possible color $c\in [k]$, let $f_c$ be the coloring of $X_{t'}$ obtained from $f$ by giving color $c$ to $v$. Also, let $H_c$ be obtained from $H$ by adding $v$, adding all monochromatic arcs incident to $v$ in $f_c$, and taking the transitive closure. The following lemma shows us how to compute entry $c_t[f,H]$.

\begin{lemma}\label{lemma:aban}
Let $G$ be a digraph, $(T,{\cal X})$ be a nice tree decomposition of $G$, and $t$ be a forget node of $T$ with child $t'$, where $X_t=X_{t'} \setminus \{v\}$. Also, let $f \in \mathcal{F}_t$ and $H\in {\cal H}_t(f)$. So, $c_t[f,H]=1$ if and only if there are $c \in [k]$ and $H'\in {\cal H}_{t'}(f_c)$ such that $c_{t'}[f_c, H']=1$, $H_c\subseteq H'$ and $H$ is an induced subgraph of $H'$.  Thus, it takes time $O(3^{tw}\cdot tw^2\cdot k)$ to decide whether $c_t[f,H]=1$.
\end{lemma}
\begin{proof}
Suppose first that $c_t[f,H]=1$ and let $g$ be an acyclic $k$-coloring of $G_t$ that extends $f$ and is represented by $H$ in $X_t$. Since $G_t=G_{t'}$, we get that $g$ is also an acyclic $k$-coloring of $G_{t'}$. Then let $c=g(v)$ and consider $H'$ the digraph representing $g$ in $X_{t'}$; by definition we know that $c_{t'}[f_c,H']=1$, where $f_c = g_{|X_{t'}}$. Again because $G_t = G_{t'}$, we get that the monochromatic paths in $G_t$ and in $G_{t'}$ are the same, which gives us that $H$ must be an induced subgraph of $H'$. Also, from Proposition~\ref{prop:trivialzero}, we get that $H'$ must contain the monochromatic arcs incident in $v$ and must be transitive, which together with the fact that $H\subseteq H'$ implies that $H_c\subseteq H'$, as we wanted to prove.

Suppose now that there are $c \in [k]$ and $H' \in {\cal H}_t(f_c)$ such that $c_{t'}[f_c,H']=1$, $H_c\subseteq H'$ and $H$ is an induced subgraph of $H'$. We want to show that $c_t[f,H]=1$, that is, we want to find an acyclic $k$-coloring $g$ of $G_t$ that extends $f$ and is represented by $H$ in $X_t$. Since $c_{t'}[f_c,H']=1$, there exists an acyclic $k$-coloring $g$ of $G_{t'}$ that extends $f_c$ and is represented by  $H'$ in $X_{t'}$. Since $G_t = G_{t'}$ and $f_c$ extends $f$, we have that $g$ is an acyclic $k$-coloring of $G_t$ that extends $f$. It remains to show that $g$ is represented by $H$ in $X_t$, that is, that for each pair of vertices $u,w \in X_t$, we get $uw \in E(H)$ if, and only if, there is a monochromatic $(u,w)$-path  by $g$ in $G_t$. Let $uw \in E(H)$. We have that $uw \in E(H')$ because $H$ is an induced subgraph of $H'$, by definition. Since $H'$ represents $g$ in $X_{t'}$, there is a monochromatic $(u,w)$-path by $g$ in $G_{t'}=G_t$, and we are done. Let now $P$ be a monochromatic $(u,w)$-path in $G_t$, with $u,w \in X_t$. Again, because $X_t\subset X_{t'}$, $G_t = G_{t'}$, and $H'$ represents $g$ in $X_{t'}$, we get that $uw\in E(H')$, and since $H$ is an induced subgraph of $H$ and $\{v,w\}\subseteq V(H)$, it follows that $uw\in E(H)$.

Finally, for each $c\in [k]$, one can compute $H_c$ in $O(tw^2)$ time, it suffices to add the monochromatic arcs to $H$ and compute the transitive closure. Then, for each supergraph of $H_c$, we check whether $H$ is an induced subgraph of $H_c$, which can be done also in $O(tw^2)$ time. By Proposition~\ref{prop:trivialzero}, we need only check acyclic supergraphs of $H_c$, and since only arcs incident to $v$ are allowed to be added to $H_c$, we have that, for each $u\in X_{t'}$ we either have $uv\in E(H')$, or $vu\in E(H')$, or neither of these arcs appear in $H'$. This gives a total of $O(3^{tw})$ supergraphs that need to be checked.
\end{proof}

Now let $t \in V(T)$ be an introduce node, with child $t'$, and let $v \in V(G)$ be the vertex that is introduced by $t$. Note that pairs of vertices of $X_{t'}$ can now be connected by a path passing through $v$. Therefore, the entry $c_t[f,H]$ might be positive due to the existence of an acyclic $k$-coloring $g$ of $G_t$ represented by $H$ in $X_t$, but that the graph $H'$ that represents $g$ in $X_{t'}$ is not necessarily equal to $H-v$. So, in order to compute $c_t[f,H]$, we need to look at additional entries of $c_{t'}$. In other words, we need to allow for some arcs of the type $uw$ where $u\in N^-_H(v)$ and $w\in N^+_H(v)$ not to exist in the investigated entries of $t'$. In what follows, we give the formal necessary definitions.

We say that $H$ is \emph{feasible} if: for every $u \in N^-_H(v) \setminus N^{-}_{G_t}(v)$, there is $w \in N^{-}_{G_t}(v)$ such that $uw \in E(H)$; and for every $u \in N^{+}_H(v) \setminus N^{+}_{G_t}(v)$, there is $w \in N^{+}_{G_t}(v)$ such that $ wu \in E(H)$. The following lemma helps us solve entries $c_t[f,H]$ that are trivially $0$.

\begin{lemma}\label{lem:feasible}
Let $G$ be a digraph, $(T,{\cal X})$ be a nice tree decomposition of $G$ and $t$ be an introduce node of $T$ with child $t'$, where $X_t=X_{t'} \cup \{v\}$. Let $f \in \mathcal{F}_t$ and $H \in \mathcal{H}_t(f)$. If $c_t[f,H] = 1$, then $H$ is feasible.
\end{lemma}
\begin{proof}
Let $g$ be an acyclic coloring of $G_t$ that extends $f$ and is represented by $H$ in $X_t$. 
Suppose that $u\in N^-_H(v)\setminus N^-_{G_t}(v)$, and let $P$ be a monochromatic $(u,v)$-path in $G_t$ (it exists since $g$ is represented by $H$). Let $w$ be the in-neighbor of $v$ in $P$, i.e., $w\in N^-_{G_t}(v)\cap V(P)$. Then $P-v$ is a monochromatic $(u,w)$-path in $G_t$, and therefore $uw\in E(H)$. The argument for the out-neighbors of $v$ in $H$ is analogous.
\end{proof}

Now, in order to outline exactly which entries from the table of $t'$ must be checked, we first define the minimum graph $H_{min}$ that can represent in $X_{t'}$ the same colorings represented by $H$ in $X_t$. Let $H_{min}$ be obtained from $H-v$ by removing all arcs from $N^-_H(v)$ to $N^+_H(v)$ that are not in $G_t$, and taking the transitive closure. Observe that, supposing that $H$ satisfies Proposition~\ref{prop:trivialzero}, we get that $H_{min}\subseteq H$.

Then, we define the family ${\cal H}^+_{min}$ of supergraphs of $H_{min}$ that are subgraphs of $H$. Note that by construction, if $H'\in {\cal H}^+_{min}$, then $E(H')\setminus E(H)$ contains only arcs between $N^-_H(v)$ to $N^+_H(v)$. 

\begin{figure}[ht]
 \begin{center}
\begin{tikzpicture}
   \tikzstyle{vertex}=[draw,thick,circle,fill=white,minimum size=7pt,inner sep=0pt]
 \hspace{-3 cm}

 \node[vertex] (1) at (0,0) {};
 \node[vertex][label=above:$v$] (v) at (1,1) {};
 \node[vertex] (2) at (0,2) {};
 \node[vertex] (3) at (2,2) {};
 \node[vertex] (4) at (2,0) {};
 

 \draw[dashed,->] (1) -- (2);
 \draw[dashed,->] (1) -- (v);
 \draw[->] (2) -- (3);
 \draw[dashed,->] (2) -- (v);
 \draw[dashed,->] (v) -- (3);
 \draw[dashed,->] (v) -- (4);

 \draw[dashed,->] (1) to[out=150,in=-135] (-0.5,2.5) to[out=45,in=105](3);
 \draw[dashed,->] (2) to[out=210,in=135] (-0.5,-0.5) to[out=315,in=255](4);

 \hspace{3cm}
 
 \node[vertex] (1) at (0,0) {};
 \node[vertex] (2) at (0,2) {};
 \node[vertex] (3) at (2,2) {};
 \node[vertex] (4) at (2,0) {};

 \draw[dashed,->] (1) -- (2);
 \draw[->] (2) -- (3);


 \hspace{3cm}

 \node[vertex] (1) at (0,0) {};
 \node[vertex] (2) at (0,2) {};
 \node[vertex] (3) at (2,2) {};
 \node[vertex] (4) at (2,0) {};

 \draw[dashed,->] (1) -- (3);
 \draw[dashed,->] (1) -- (2);
 \draw[->] (2) -- (3);

\end{tikzpicture}

\vspace{.5cm}

 \caption{To the left we have the graph $H$ and all dashed arcs are not in $G_t$, we also only depict the component of $H$ containing $v$, because it is the only relevant one to the construction of $H_{min}$. In the middle, a partial graph $H'$ that consists of arcs that are not from $N^-_H(v)$ to $N^+_H(v)$ along with the chromatic arc in $G_t$. To the right, we have $H_{min}$ that is obtained from $H'$ by taking its transitive closure.}
\label{fig:Hmin}
 \end{center}
\end{figure}
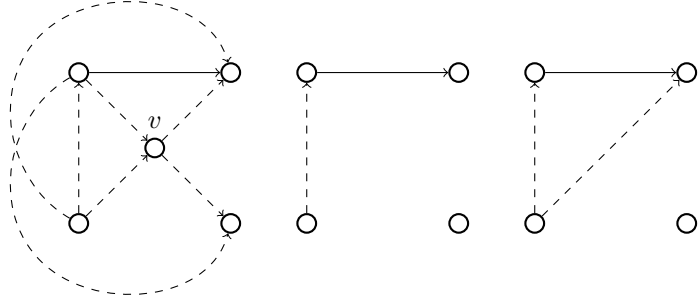

An example of the definition of $H_{min}$ can be seen in Figure \ref{fig:Hmin}.

Given a function $f: A \rightarrow C$, defined on a set $A$ and let $B$ be a subset of $A$, we define the restriction of $f$ to $B$ and denote it by $f_{|_{B}}$ as the function $f_{|_{B}}:B \rightarrow C$ such that $f_{|_{B}}(x) = f(x)$ for every $x \in B$. Also, if $ g: B \rightarrow C$ is the restriction to $B$ of some function $f: A \rightarrow C$, we say that $f$ extends $g$ or that $f$ is an extension of $g$.

\begin{lemma}\label{lemma:intro}
Let $G$ be a digraph, $(T,{\cal X})$ be a nice tree decomposition of $G$, and $t$ be an introduce node of $T$ with child $t'$, where $X_t=X_{t'} \cup \{v\}$. Also, let $f \in \mathcal{F}_t$ and $H\in {\cal H}_t(f)$. So $c_t[f,H]=1$ if and only if $H$ is feasible, acyclic, transitive, $H_f\subseteq H$, and  there is $H' \in \mathcal{H}^+_{min}$ such that $c_{t'}[f',H'] =1$, where $f'=f_{|{X_{t'}}}$. Thus, it takes time $O(2^{tw}\cdot tw^2)$ to decide whether $c_t[f,H]=1$.

\end{lemma}

\begin{proof}
Suppose first $c_t[f,H]=1$. Proposition~\ref{prop:trivialzero} and Lemma~\ref{lem:feasible} give us that $H$ is feasible, acyclic, transitive, and $H_f\subseteq H$.  
Now, let $g$ be an acyclic $k$-coloring of $G_t$ that extends $f$ and is represented by $H$ in $X_t$. Let $f'=g_{|{X_{t'}}}$ (which is also equal to $f_{|X_{t'}}$), and $g' = g_{|V(G_{t'})}$. Since $G_{t'} = G_t-v$, we know that $g'$ is an acyclic $k$-coloring of $G_{t'}$ that extends $f'$, so it remains to show that $g'$ is represented by a digraph in ${\cal H}^+_{min}$. Let $H'$ be the digraph that represents $g'$ in $X_{t'}$. 
Because $G_{t'} = G_t-v$, we know that every path in $G_{t'}$ also exists in $G_t$, and hence $H'\subseteq H$; thus we need only to prove that $H_{min}\subseteq H'$.

Since $V(H_{min})=V(H')=X_{t'}$, we need to show that $E(H_{min})\subseteq E(H')$. So, let $uw\in E(H_{min})$. 
Recall that $H_{min}$ is the transitive closure of $H'_{min}$, where $H'_{min}$ is obtained from $H-v$ by removing all arcs between $N^-_H(v)$ and $N^+_H(v)$ that are not in $G_t$. This gives us that $uw\notin E(G_t)$ since $H'$ contains all monochromatic arcs of $G[X_{t'}]$. 
Thus $uw$ was introduced when taking the transitive closure of $H'_{min}$, which implies that either there exists $z\in X_{t'}$ such that $uz,zw\in E(H'_{min})$, or there exist $u'\in N^-_{H}(v)$ and $w'\in N^+_H(v)$ such that $uu',w'w\in E(H'_{min})$ and $u'w'\in E(G_t)$. In either case, observe that all paths represented by arcs of $H$ that are not between $N^-_H(v)$ and $N^+_H(v)$ still exist in $G_{t'}$. This means that we can combine either the $(u,z)$-path and the $(z,w)$-path into a $(u,w)$-path in $G_{t'}$ in the former case, or, in the latter case, the $(u,u')$-path, the $(w',w)$-path and the arc $u'w'$ into a $(u,w)$-path in $G_{t'}$. It thus follows that $uw\in E(H')$, as we wanted to prove.

Suppose now that $H$ is feasible, acyclic, transitive, $H_f\subseteq H$, and there is $H' \in \mathcal{H}^+_{min}$ such that $c_{t'}[f',H']=1$, where $f'=f_{|{X_{t'}}}$. By the definition of $c_{t'}[f',H']$, there is an acyclic $k$-coloring $g$ of $G_{t'}$ that extends $f'$ and is represented by $H'$. Let $g$ be the $k$-coloring of $G$ obtained from $g'$ by coloring $v$ with $f(v)$. 
Since $g$ extends $f$, we need to prove that $g$ is an acyclic $k$-coloring of $G_t$ represented by $H$.
First suppose by contradiction that $g$ produces a monochromatic cycle $C=(v_1, v_2, \ldots, v_q)$. We must have $v \in V(C)$, otherwise $C$ would be contained in $G_{t'}$, contradicting the fact that $g'$ is an acyclic $k$-coloring of $G_{t'}$. Suppose then, without loss of generality, that $v=v_1$. Since $N_{G_t}(v)\subseteq X_t$, we have that $v_2,v_q \in X_t$; therefore $vv_2,v_q v $ are monochromatic edges of $G[X_t]$, and since $H_f\subseteq H$ and $H$ is transitive, we get $\{v_qv,vv_2,v_qv_2\} \subseteq E(H)$. Furthermore, because the path $(v_2,\ldots, v_q)$ is in $G_{t'}$ and $H'$ represents $g'$ in $X_{t'}$, we have that $v_2v_q \in E(H')$. But $H'\subseteq H$ by the definition of ${\cal H}^+_{min}$, which contradicts the fact that $H$ is acyclic as in this case we have the cycle $(v_2,v_q,v_2)$ in $H$. 
Therefore, $g$ is an acyclic $k$-coloring of $G_t$. 

Now we show that $H$ represents $g$, that is, that $uw \in E(H)$ if and only if there is a monochromatic $(u ,w)$-path in $G_t$. 
So first consider $uw \in E(H)$. If $uw\in E(H')$, then there is nothing to do as $H'$ represents $g'$ in $X_{t'}$. So suppose $uw \in E(H) \setminus E(H')$. We analyze the following cases:
    
    \begin{enumerate}
    
        \item[$(a)$] $w=v$: if $u \in N^-_{G_t}(v)$, then $(u,v)$ is the path we are looking for, so suppose $u \notin N^-_{G_t}(v)$. Since $H$ is feasible, there exists $z \in N^-_{G_t}(v)$ such that $uz \in E(H)$. As $\{u,z\} \subseteq N^-_H(v)$, we have that $uz \in E(H')$. So there is a $(u,z)$-path $P$ in $G_{t'}$ which concatenated with the arc $zv$ gives us a $(u,v)$-path in $G_t$, as desired.
        
        \item[$(b)$] $u=v$: analogous to item $(a)$.
        
        \item[$(c)$] $v \notin \{u,w\}$: by construction, we know that $u \in N^-_H(v) $ and $w \in N^+_H(v)$. 
        Applying items $(a)$ and $(b)$, we get a $(u,v)$-path and a $(v,w)$-path which concatenated yield the desired $(u,w)$-path.
        
    \end{enumerate}

Now, let $P$ be a monochromatic $(u,w)$-path in $G_t$. If $v \notin V(P)$, then we are done since $H'$ represents $g'$ in $X_{t'}$, which implies that $uw \in E(H')$, and because $H'\subseteq H$ by the definition of ${\cal H}^+_{min}$. So suppose $v\in V(P)$ and write $P=(u=v_1, \ldots , v=v_i, \ldots , v_q = w)$. Since $(v_1,\ldots, v_{i-1})$ and $(v_{i+1},\ldots, v_{q})$ are paths in $G_{t'}$ and $H'$ represents $g'$ in $X_{t'}$, we have that $\{v_1v_{i-1},v_{i+1}v_q\} \subseteq E(H')\subseteq E(H)$. Also, since $v_{i-1} \in N^-_{G_t}(v)$, $v_{i+1} \in N^+_{G_t}(v)$, and $H_f\subseteq H$, we get that $\{v_{i-1}v,vv_{i+1}\} \subseteq E(H)$. It follows that $v_1v_q \in E(H)$ as $H$ is transitive.

Finally, checking whether $H$ is feasible, acyclic, transitive and $H_f\subseteq H$ can be done in $O(tw^2)$ time. Also, computing $H_{min}$ can be done in $O(tw^2)$. Observe that, because $H$ is acyclic, we have that $N^-_H(v)$ and $N^+_H(v)$ are disjoint. Hence, since $E(H')\setminus E(H)$ contains only arcs from $N^-_H(v)$ to $N^+_H(v)$, we get that ${\cal H}^+_{min}$ has at most $2^{tw}$ digraphs: for each pair $u\in N^-_H(v)$ and $w\in N^+_H(v)$, either arc $uw$ is in $H'\in {\cal H}^+_{min}$, or it is not. 
\end{proof}

It remains to investigate join nodes, and for this we need the following definition. Let $t$ be a join node, with children $t_1$ and $t_2$. Given $H_1 \in \mathcal{H}_{t_1}(f)$ and $H_2\in \mathcal{H}_{t_2}(f)$, we say that $H_1$ and $H_2$ \emph{combine into $H$} if $H$ is the transitive closure of $H_1 \cup H_2$, where $H_1\cup H_2$ is simply the digraph $(X_t, E(H_1)\cup E(H_2))$ (recall that $X_{t_1} = X_{t_2} = X_t$). 

\begin{lemma}\label{lemma:jun}
Let $G$ be a digraph, $(T,{\cal X})$ be a nice tree decomposition of $G$ and $t$ be a join node with children $t_1,t_2 $. Let $f \in \mathcal{F}_t$ and $H\in {\cal H}_t(f)$. So $c_t[f,H] =1$ if and only if $H$ is acyclic and transitive and there are $H_1 \in \mathcal{H}_{t_1}(f)$ and $H_2 \in \mathcal{H}_{t_2}( f)$ that combine into $H$ such that $c_{t_1}[f,H_1] = c_{t_2}[f,H_2] =1$. Thus, it takes time $O(9^{tw}\cdot tw^2)$ to decide whether $c_t[f,H]=1$.
\end{lemma}

\begin{proof}
Suppose first that $c_t[f,H]=1$, and let $g$ be an acyclic $k$-coloring of $G_t$ that extends $f$ and is represented by $H$ in $X_t$. Since for each $i \in [2]$, $G_{t_i}$ is a subgraph of $G_t$, we have that $g_i$ equal to $g$ restricted to $G_{t_i}$ is an acyclic $k$-coloring of $G_{t_i}$. Let $H_i$ be the digraph that represents $f_i$ in $X_{t_i}=X_t$; then, $c_{t_i}[f,H_i]=1$, for each $i \in [2]$. 
It remains to show that $H_1,H_2$ combine into $H$. For this, let $H'$ be the transitive closure of $H_1 \cup H_2$; we prove that $H = H'$. 
First consider $uv \in E(H)$, and let $P$ be a monochromatic $(u,v)$-path in $G_t$ (it exists since $H$ represents $g$ in $X_t$). 
If $P \subseteq G_{t_i}$ for $i \in [2]$, we have that $uv \in E(H_i)$ and therefore $uv \in E(H')$. So suppose that $P$ is not contained in $G_{t_i}$, for each $ i \in [2]$. 
Write $P$ as $(v_1=u,\ldots,v_p=v)$, and note that there must exist $i_1 < \ldots < i_q$, $q\ge 1$, such that $\{v_{i_1},\ldots,v_{i_p}\}\subseteq V(P)\cap X_t$ and these vertices break $P$ into pieces that alternate among $G_1$ and $G_2$; more formally, we have that, without loss of generality, $P_1 = (v_1,\ldots,v_{i_1})\subseteq G_1$, $P_j = (v_{i_j},\ldots,v_{i_{j+1}})\subseteq G_1$ for every even $j\in [q]$, and $P_j = (v_{i_j},\ldots,v_{i_{j+1}})\subseteq G_2$ for every odd $j\in [q]$, with $i_{q+1}$ being equal to $q$. Since $H_i$ represents $g_i$ in $X_{t_i}$, for each $i\in [2]$, this gives us that the path $(v_1,v_{i_1},\ldots,v_{i_q},v_p)$ is contained in $H_1\cup H_2$, which in turn gives us that $v_1v_q = uv\in E(H')$. 
Now let $uv \in E(H')$. If $uv \in E(H_i)$, for some $i \in [2]$, then there is a monochromatic $(u,v)$-path in $G_{t_i} \subseteq G_t$. Since $H$ represents $g$ in $X_t$, we have $uv \in E(H)$. So, suppose that $uv\in E(H')\setminus E(H_1\cup H_2)$, i.e., $uv$ is added when taking the transitive closure of $H_1 \cup H_2$. So, there is a $(u,v)$-path $(u=w_1, \ldots , w_q =v )$ in $H_1 \cup H_2$. Since $H_i$ represents $g_i$ in $X_{t_i}$ for each $i \in [2]$, we have a $(w_j,w_{j+1})$-path $P_j$ contained in either $G_{t_1}$ or in $G_{t_2}$, for each $j\in [q-1]$. The concatenation of these paths thus contain a monochromatic $(u,v)$-path in $G_t$, and hence $uv \in E(H)$.

Now suppose that there are $H_1 \in \mathcal{H}_{t_1}(f)$ and $H_2 \in \mathcal{H}_{t_2}(f)$ that combine into $H$ and such that $c_{t_1}[f,H_1 ] = c_{t_2}[f,H_2] =1$. 
For each $i \in [2]$, let $g_i$ be an acyclic $k$-coloring of $G_{t_i}$ that extends $f$ and is represented by $H_i$ in $X_{t_i} = X_t$. We prove that $g$ equal to the union of $g_1$ and $g_2$ (i.e., $g(v) = g_i(v)$ for each $v\in V(G_i)$, for each $i\in [2]$) is an acyclic $k$-coloring of $G_t$ represented by $H$. 
To see that $g$ is acyclic, one can suppose by contradiction that there is a monochromatic cycle $C$ in $G_t$ and apply arguments similar as before to conclude that each part of $C$ contained in $G_i$ defines an edge in $H_i$. Combining these edges we end up with a cycle in $H_1\cup H_2\subseteq H$, contradicting the fact that $H$ is acyclic. 
The proof that $g$ is represented by $H$ is much similar to the proof in the previous paragraph, so we refrain from repeating it.

Finally, for each pair $H_1\in {\cal H}_{t_1}(f)$ and $H_2\in {\cal H}_{t_2}(f)$, testing whether $H_1$ and $H_2$ combine into $H$ can be done in $O(tw^2)$, it suffices to take the transitive closure of $H_1\cup H_2$ and test whether it is equal to $H$. Since we are interested only in pairs for which the corresponding entries are~1, by Proposition~\ref{prop:trivialzero} we can constrain ourselves to acyclic digraphs $H_1$, $H_2$ whose vertex sets is $X_t$. There are at most $3^{tw}$ such digraphs, as for each $u,v\in X_t$, either $uv$ is an arc, or $vu$ is an arc, or neither of them is. Hence, there is a total of $(3^{tw})^2 = 9^{tw}$ such pairs.
\end{proof}

We now combine the previous lemmas to produce the main result of this section.

\begin{theorem}\label{thm:alg}
Let $G$ be a digraph and $k$ a positive integer. If a nice tree decomposition of the underlying graph of $G$ with width $tw$ is known, then one can decide whether $G$ has an acyclic $k$-coloring  in time $\mathcal{O}(k^{tw}\cdot 27^{tw^2}\cdot tw^2\cdot n)$.
\end{theorem}

\begin{proof}
As is standard, the computation is made in a bottom-up manner. 
The time spent on a node $t$ is given by the size of the table $c_t$, multiplied by the time spent on each entry of $c_t$. From previous lemmas, we know that computing an entry takes time at most $O(9^{tw}\cdot tw^2)$. For the size of the tables, first note that $|{\cal F}_t| = k^{tw}$ since there are $k$ possibilities for each $u\in X_t$. As for $|{\cal H}_t(f)|$, we can consider only those digraphs that are acyclic, due to Proposition~\ref{prop:trivialzero}. Thus, for each pair $u,v\in X_t$, we have 3 possibilities: either only the arc $uv$ is in $H$; or only the arc $vu$ is in $H$; or neither of them is in $H$. We therefore have a total of at most $3^{tw^2}$ digraphs to be considered, giving us a table $c_t$  of size $k^{tw}\cdot 3^{tw^2}$. Because there are $O(k\cdot n)$ nodes in a nice tree decomposition~\cite{CMF.15}, we get the desired time complexity.
\end{proof}

    It is known that a tree decomposition of a graph $G$ with a treewidth at most $w=5tw(G)+4$ can be computed in time~$2^{\mathcal{O}(w)}\cdot n$~\cite{B.16}. It is also known that $\chi_a(G)\le \chi(G)\le tw(G)+1$~(folklore). Hence, in order to compute $\chi_a(G)$, we can apply Theorem \ref{thm:alg} for each value $k\in [w]$, which would take time $\mathcal{O}(w^{w+3}\cdot 27^{w}\cdot n)$. The corollary below thus follows.

\begin{corollary}
    Let $G$ be any digraph. One can compute the dichromatic number of $G$ in {\FPT} time when parameterized by the treewidth of $G$.
\end{corollary}


\section{Final Remarks}\label{sec:conclusion}

In this paper, we have proved that a number of results on the chromatic number of graphs can be generalized for the dichromatic number of digraphs. Additionally, we have investigated the dichromatic number of products of dicycles, given exact values for almost all possible combinations, as discussed shortly. Finally, we have also given exact values for the chromatic number of the odd product of a cycle by a bipartite graph, thus complementing previous results, as discussed in Section~\ref{sec:prelim}. There is still a lot to be learnt, even about the chromatic number of products. In what follows, we propose a few questions that we found more interesting, given the results presented here.

As we have discussed in Section~\ref{sec:prelim}, the famous Hedetniemi's Conjecture has been recently disproved~\cite{S.19}. 
Because acyclic colorings of digraphs are generalizations of proper colorings of  graphs, it also follows that equality is not valid for the dichromatic number. However, the counterexamples obtained by~\cite {S.19}, being graphs, imply that the counterexamples for equality on digraphs have dicycles of length~2. An interesting question, therefore, is whether equality holds in the case of oriented graphs. 

\begin{question}
Let $G$ and $H$ be oriented graphs. Is it true that $\cn(G\times H) = \min\{\cn(G),\cn(H)\}$?
\end{question}

Observe that, by Theorem~\ref{thm:positiveAnswer}, we get a positive answer for products of dicycles. In fact, we get the following corollary.

\begin{corollary}\label{cor:directCycles}
Let $\vv{C}_{n_1}, \ldots, \vv{C}_{n_t}$ be dicycles on $n_1,\ldots , n_t$ vertices, respectively. Then, $$\cn(\vv{C}_{n_1} \times \ldots \times \vv{C}_{n_t})=2.$$
\end{corollary}

Still concerning dicycles, one can observe that exact values are known for each of the products when both factors are dicycles (see Corollaries~\ref{cor:DAG} and~\ref{cor:directCycles}, and Theorems~\ref{thm:Cartesian},~\ref{thm:prodforte} and~\ref{thm:nlechi}). 
Additionally, it is easy to determine $\cn(\vv{C}_n\square H)$ and $\cn(\vv{C}_n\times H)$ in terms of $\cn(H)$. Indeed, if $H$ is acyclic, then $\cn(\vv{C}_n\square H) = 2$ and $\cn(\vv{C}_n\times H) = 1$; otherwise $\cn(\vv{C}_n\square H) = \cn(H)$ and $\cn(\vv{C}_n\times H) = 2$. Therefore, when only one of the factors is a dicycle, the remaining open cases are stated in the following question.

\begin{question}
Let $H$ be any digraph and $n$ be a positive integer. Can $\cn(H\boxtimes \vv{C}_n)$ and/or $\cn(H[\vv{C}_n])$ be written as a function of $\cn(H)$ and $n$?
\end{question}

Finally, observe that the values obtained in Theorem~\ref{thm:prodforte} are exactly equal to the upper bound given by $\cn(\vv{C}_n[H])$ in Theorem~\ref{thm:nlechi}. We then ask whether this upper bound is tight in general.

\begin{question}
Let $H$ be any digraph and $n$ be a positive integer. Does it hold that $\cn(\vv{C}_n\boxtimes H) = \cn(H)+ \left\lceil\frac{\cn(H)}{n-1}\right\rceil$?
\end{question}

\bibliography{main}



\appendix

\smallskip

 \end{document}